\documentclass[a4paper,11pt,english]{smfart}
\usepackage{amsfonts}
\usepackage{amsmath,amsthm}
\usepackage{latexsym}
\usepackage{array}
\usepackage{amssymb}
\usepackage{graphicx}
\usepackage{enumerate}
\usepackage[english]{babel}

\usepackage{color}
\definecolor{marin}{rgb}   {0.,   0.3,   0.7} 
\definecolor{rouge}{rgb}   {0.8,   0.,   0.} 
\definecolor{sepia}{rgb}   {0.8,   0.5,   0.} 
\usepackage[colorlinks,citecolor=marin,linkcolor=rouge,
            bookmarksopen,
            bookmarksnumbered
           ]{hyperref}

\newtheorem{lemma}{Lemma}[section]

\newtheorem{remark}[lemma]{Remark}
\newtheorem{example}[lemma]{Example}

\newtheorem{notation}[lemma]{Notation}
\newtheorem{definition}[lemma]{Definition}
\newtheorem{conclusion}[lemma]{Conclusion}
\numberwithin{equation}{section}
\newcommand{\QED}{\mbox{}\hfill \raisebox{-0.2pt}{\rule{5.6pt}{6pt}\rule{0pt}{0pt}} 
          \medskip\par}

\newcommand{\dd}{\,\mathrm{d}}

\newcommand{\Hc}{\mathcal{H}}

\newcommand{\R}{\mathbb{R}}

\newcommand{\Z}{\mathbb{Z}}

\let\oldmarginpar\marginpar
\renewcommand\marginpar[1]{\-\oldmarginpar[\raggedleft\footnotesize #1]%
{\raggedright\footnotesize #1}}

\author{Nicolas Crouseilles}
\address{INRIA \& IRMAR  \\
 F-35042 Rennes, France. } 
\email{nicolas.crouseilles@inria.fr}

\author{Lukas Einkemmer}
\address{University of Innsbruck  \\
Technikerstra\ss e 19a, 
A-6020 Innsbruck, Austria. } 
\email{lukas.einkemmer@uibk.ac.at}

\author{Erwan Faou}
\address{INRIA \& ENS Cachan Bretagne  \\
Avenue Robert Schumann F-35170 Bruz, France. } 
\email{Erwan.Faou@inria.fr}

\title[AP scheme for the relativistic Vlasov--Maxwell equations]
{An asymptotic preserving scheme for the relativistic Vlasov--Maxwell equations in the classical limit}

\begin{document}

\begin{abstract}
We consider the relativistic Vlasov--Maxwell (RVM) equations in the limit when the light velocity $c$ goes to infinity.
In this regime, the RVM system converges towards the Vlasov--Poisson system and the aim of this paper is to construct asymptotic preserving numerical schemes that are robust with respect to this limit. 

Our approach relies on a time splitting approach for the RVM system employing an implicit time integrator for Maxwell's equations in order to damp the higher and higher frequencies present in the numerical solution. 

A number of numerical simulations are conducted in order to investigate the performances of our numerical scheme both in the relativistic as well as in the classical limit regime. In addition, we derive the dispersion relation of the Weibel instability for the continuous and the discretized problem.
\end{abstract}

\thanks{This work is supported by the Fonds zur F\"orderung der Wissenschaften (FWF) -- project id: P25346.
}
\thanks{This work was partly supported by the ERC Starting Grant Project GEOPARDI No 279389 and by the ANR project Moonrise ANR-14-CE23-0007-01. }
\thanks{The computational results presented have been achieved in part using the Vienna Scientific Cluster (VSC)}
\maketitle
\tableofcontents

\section{Introduction}

\textcolor{black}{In a wide range of plasma processes, interactions between relativistic particles and electromagnetic fields play a very important role. For instance, it is possible to accelerate particles to relativistic speed (see \cite{laser} for a discussion of recent experiments in laser-plasma interaction). In this context, the full relativistic Vlasov--Maxwell (RVM) system is indispensable for dealing with the deviation from thermal equilibrium as well as to describe magnetic effects. Nevertheless, if the speed of light is large compared to the characteristic thermal speed of the plasma (a situation often encountered in practice, see \cite{apart, drouin, duclous}), standard numerical methods have to take so 
small time steps as to render them extremely inefficient when applied to such problems. Moreover, using the Vlasov--Poisson equation is not an option in this case as this model completely neglects all magnetic effects. Certainly a code that performs well in the situation described would be useful. This enables the study of the range of parameter values for which relativistic effects (such as encountered in the Weibel instability) are essential and the range for which electrostatic effects (such as Landau damping) take over. }

In this work, we are interested in non relativistic limit of the relativistic Vlasov--Maxwell (RVM) equations. This system of nonlinear partial differential equations couples a Maxwell system with a transport equation for the particles density and depends on a parameter $c$ which represents the speed of light. It has been shown in \cite{schaeffer_cmp86, asano, degond} that for smooth initial data  with compact support  classical solutions exist on an intervals $[0, T]$ independent of $c$  and converge  
to the solution of the Vlasov--Poisson system at a rate proportional to $1/c$ as $c$ tends to infinity. \textcolor{black}{Note that we exclusively consider this so-called electric limit in the present work (for more details on the so-called 
electric and the magnetic limits, we refer the reader to \cite{manfredi2013}).}

To reproduce this behavior numerically, standard schemes usually require very small time steps since solutions develop highly oscillatory phenomena on a scale proportional to $1/c$. 
The main goal of this work is to overcome this difficulty by deriving numerical schemes that preserve this Vlasov-Poisson limit without requiring small time steps. 

Starting with the seminal paper of Cheng \& Knorr \cite{cheng}, a largy body of works has been devoted to the solution of the Vlasov--Poisson system (see, for example, \cite{chengy, ccfm, einkemmerdg2014, heath2012}). Recently, both the relativistic as well as the non-relativistic Vlasov--Maxwell system has received some attention (see, for example, \cite{cglm, cef, cns, mangeney, valis, suzuki, besse_latu, shebalin}). 
Regarding time integration, splitting methods have several advantages: They are often explicit and computationally attractive as they reduce the integration of the system to a sequence of numerical approximations of lower dimensional problems, and in general structure preserving (symplecticity, reversibility, see \cite{blanes2008,hlw} for general settings). 
 For example in the Vlasov--Poisson case the computational advantage lies in the fact that splitting methods  reduce the nonlinear system to a sequence of one-dimensional explicit advections.
Various space discretization methods can then be employed to solve the resulting advections.
 Semi-Lagrangian methods using interpolation with Fourier or spline basis functions as well as discontinuous Galerkin type schemes are among the most commonly employed methods.

 In our previous work \cite{cef} we have introduced a three-term splitting for the Vlasov--Maxwell system that is computationally attractive, easy to implement, and extensible to arbitrary order in time. In addition, it can be easily combined with a range of space discretization techniques. \textcolor{black}{This numerical scheme is based on the Poisson bracket proposed in \cite{Morrison80} which, as has been pointed out in \cite{qin2015}, does not satisfy the Jacobi identity, see \cite{MaWe82}. Nevertheless, the numerical method introduced in \cite{cef} is time reversible, and preserves the Poisson equation as well as the divergence free condition on the magnetic field. In addition, the numerical results given in \cite{cef} show that it is superior with respect to energy conservation and shows better qualitative results compared to other methods from the literature. }

 \textcolor{black}{The method introduced in \cite{cef} can be easily extended to the fully relativistic case.}
However, since the scheme relies on an explicit time stepping scheme for Maxwell's equations, 
significant difficulties appear when $c$ is large (the CFL condition is proportional to $1/c^2$). 

Our goal is to design a numerical scheme that is uniformly efficient both when $c$ is of order one and 
for arbitrary large values of $c$, with a fixed set of numerical parameters. This is the context of {\em asymptotic preserving 
schemes} (see \cite{jin}). The main idea is to propose a modification of the splitting introduced in \cite{cef} 
to capture the correct asymptotic behavior without destroying the order of convergence on the limit system. 
More precisely, the linear part of Maxwell's equations 
is solved by using an implicit numerical scheme (implicit Euler or the Radau IIA method). While this choice destroys the geometric structure of the splitting (reversibility and symplecticity), the resulting time integrator is 
unconditionally stable with respect to $c$ and introduces enough numerical dissipation 
to recover the correct limit of the system when $c$ tends to infinity. This way, our scheme enjoys the 
asymptotic preserving property. 

Regarding space approximation, we will mainly consider an approach based on Fourier techniques. Let us emphasize, however, that our numerical scheme could be easily extended to various other space discretization methods. 

In section \ref{sec:vm} we will discuss the Vlasov--Maxwell system as well as its asymptotic behavior. The numerical method proposed in this paper is introduced in section \ref{sec:method}. In section \ref{sec:numerical-results} we present the numerical simulations used to benchmark and validate our asymptotic preserving scheme. Finally, we conclude in section \ref{sec:conclusion}.

\section{Relativistic Vlasov--Maxwell system \& Asymptotic behavior \label{sec:vm}}

We consider the Vlasov--Maxwell system that is satisfied by a electron distribution function  $f=f(t, x, v)$ and the 
electromagnetic fields $(E,B)=(E(t, x), B(t, x))\in\mathbb{R}^3\times \mathbb{R}^3$. 
Here, the spatial variable is denoted by  $x\in X^3$ ($X^3$ being a three dimensional torus),  
 the velocity/momentum variable is denoted by $p\in\mathbb{R}^3$, and the time is denoted by $t\geq 0$. 
Using dimensionless units, the Vlasov--Maxwell system can be written as 
\begin{equation}
\label{eq:vm}
\begin{array}{l}
\displaystyle \partial_t f + \frac{p}{\gamma} \cdot \nabla_x f + \left(E  + \frac{p}{\gamma} \times B \right) \cdot \nabla_p f = 0, \\
\partial_t E = c^2 \displaystyle \nabla_x \times B  - \int_{\R^3} \frac{p}{\gamma} f(t,x,p) \dd p + {\bar{J}(t)}, \\
\partial_t B = - \nabla_x \times E,
\end{array}
\end{equation}
where
\begin{equation}
\label{jbar}
\bar{J}(t) = \frac{1}{|X^3|} \int_{X^3} \int_{\R^3}  \frac{p}{\gamma} f(t,x,p) \dd x \dd p
\end{equation}
and $|X^3|$ denotes the volume of $X^3$. In the relativistic case, the Lorentz factor depends on $p$ and the dimensionless parameter $c$,  and is given by 
$$
\gamma=\sqrt{1+|p|^2/c^2}.
$$ 
Let us note that the splitting method considered in \cite{cef} applies to the case $\gamma=1$ and   
 $c=1$, but as we will see later, it can be easily extended to the 
case $\gamma=\sqrt{1+|p|^2/c^2}$.

In addition, two constraints on the electromagnetic field $(E, B)$ are imposed
 \begin{equation}
\label{eq:vm2}
\nabla_x\cdot  E = \rho(t,x) := \color{black} \int_{\R^3} f(t, x, p) dp -1, \;\; \nabla_x \cdot B=0, 
\end{equation}
and we easily check that if these constraints are satisfied at the initial time, they are satisfied for all times $t>0$.
We moreover impose that $E$ and $B$ are of zero average for all times $t$; that is 
\begin{equation}
\label{tgv}
\int_{X^3} E(t,x) \dd x = \int_{X^3} B(t,x) \dd x = 0, 
\end{equation}
which implies the presence of $\bar{J}$ in the system \eqref{eq:vm}. 
Moreover, the total mass is preserved; that is the relation  
\[ \int_{X^3}\int_{\R^3}f(t, x, p)\dd x \dd p = \int_{X^3}\int_{\R^3}f(0, x, p)\dd x \dd p=1\]
holds true for all times $t>0$. Let us duly note, however, that the constraints considered above are not always satisfied for a given numerical approximation. 

As initial condition, we have to specify the distribution function and the field variables: 
\[ f(t=0, x, v)=f_0(x, v), \quad B(t=0, x)=B_0(x), \]
where $\nabla_x\cdot B_0(x)=0$  and  $E(t=0,x)$ is determined by solving the Poisson equation at $t=0$ (see \eqref{eq:vm2}). 

The Hamiltonian associated with the Vlasov--Maxwell system is given by (see \cite{Morrison80, MaWe82})
\begin{eqnarray}
\Hc &:=& \frac12 \int_{X^3} |E|^2 \dd x + \frac{c^2}{2} \int_{X^3} |B|^2 \dd x + c^2 \int_{X^3 \times \R^3} \left[ \gamma -1\right] f \dd p \dd x\nonumber\\
\label{ham-split}
&=:& \Hc_E + \Hc_B + \Hc_f. 
\end{eqnarray}
The three terms correspond to electric energy, magnetic energy, and kinetic energy, respectively.  We easily check that this total energy is conserved along the exact solution of \eqref{eq:vm}.

In the limit $c\rightarrow +\infty$, the Vlasov--Maxwell equations lead to  the 3-dimensional Vlasov--Poisson 
equations (see \cite{degond, schaeffer_cmp86, asano, bostan}). Formally,  when $c$ goes 
to infinity, we check from \eqref{eq:vm} that $B$ goes to zero (assuming a well-prepared initial condition; for example,  
$B(t=0, x) = \mathcal{O}(1/c)$). In addition, $\gamma$ converges to $1$, so that 
we obtain the so-called Vlasov--Amp\`ere model %
\begin{equation}
\label{eq:vp}
\partial_t f + p\cdot \nabla_x f + E\cdot \nabla_p f = 0,  \;\; \partial_t E = - (J-\bar{J}), 
\end{equation}
with $J=\int_{\R^3} p f \dd p$. Since in the limit the electric field is curl free (i.e.,~$\nabla_x \times E=0$), we deduce that there exists a potential $\phi$ such that 
 $E=-\nabla_x \phi$. Assuming that the Poisson equation $-\Delta_x \phi = \int_{\R^3} f \dd p -1$ is satisfied for $t=0$,  we can use the continuity equation (obtained by integrating \eqref{eq:vp} with respect to $p\in \R^3$) 
 $$
 \partial_t \rho +\nabla_x\cdot J = 0, 
 $$
to verify that the Poisson equation holds true for all times $t>0$, even if we only assume that $E$ satisfies the Amp\`ere 
 equation $\partial_t E=-(J-\bar{J})$.  As a consequence, the Vlasov--Amp\`ere equation \eqref{eq:vp} is equivalent to the Vlasov--Poisson model. 
 
 \textcolor{black}{
Let us remark that the "semi-relativistic" case which corresponds to setting $\gamma=1$ in \eqref{eq:vm} 
will be also considered in this work. The limit $c\rightarrow +\infty$, i.e.~that the speed of light is large compared to 
the characteristics velocities of the problem, yields the Vlasov-Poisson equations (similar to the fully relativistic case). }

Let us note that Maxwell's equations support plane wave solutions of the form
$e^{i(k \cdot x-\omega t)}$
for $\omega=c|k|$, $k$ being the dual Fourier variable on the torus $X^3$. This is called the dispersion relation. We thus conclude that, for a fixed wavenumber $k$, the angular frequency $\omega$ increases proportional to $c$. This poses a difficult problem for a given numerical scheme as high frequency oscillations have to be resolved. This is especially important as the nonlinear coupling to the Vlasov equation excites modes that are not present in the initial condition. Of course, in the latter case the dispersion relation is modified as we have to take the full Vlasov--Maxwell system into account. This point will be further discussed in section \ref{sec:dispersion}.

 \section{Description of the numerical method \label{sec:method}}
 In this section, we propose a time discretization of \eqref{eq:vm} enjoying the asymptotic preserving 
 property in the sense that it is uniformly stable with respect to $c$ and is consistent with the Vlasov--Poisson model \eqref{eq:vp} when $c$ goes to infinity, for a fixed time step. 
 We first focus on the time discretization of the linear part of the Maxwell's equations before describing the 
 time discretization of the rest of the RVM model. Then, a fully discretized presentation of the numerical scheme 
 is performed in the case of the 1+1/2 RVM model. 
 
 In the sequel, we will use a discretization of the time variable $t^n=n\Delta t$, $\Delta t>0$ and the classical 
notation $Z^n$ as approximation of $Z(t^n)$ where $Z$ can denote the electric (or magnetic)  field $E$ ($B$) 
 as well as the distribution function $f$. Finally, in the third part of this section we will denote
 the Fourier transform of any space dependent quantity $Z$ by $\hat{Z}$ and the associated frequency in Fourier space by $k$. 
 
 \subsection{Time discretization of Maxwell's equations}

We split Maxwell's equations between the linear part 
\begin{equation}
\label{linear_maxwell}
 \partial_t E = c^2\nabla_x\times B, \qquad \partial_t  B=-\nabla_x\times E, 
\end{equation}
 and the nonlinear part 
 $$
 \partial_t E=-(J-\bar{J}). 
 $$
 The former is essentially a wave equation, which is stiff due to the presence of $c$), 
 while the latter is a (non stiff) nonlinear part that only mediates the coupling to the Vlasov equation 
 and will be considered in the next section. 

 In order to avoid the stability constraint imposed by $c^2$, an implicit scheme has to be used for the linear part \eqref{linear_maxwell}. 
 Let us consider an implicit Euler scheme  
 \begin{equation}
\label{maxwell-linear-time}
 E^{n+1}=E^n + c^2\Delta t \nabla_x\times B^{n+1}, \quad  B^{n+1}=B^n-\Delta t\nabla_x\times E^{n+1}, 
\end{equation}
 such that combining the two equations gives an implicit time discretization for the wave equation satisfied by $B$ 
 \begin{eqnarray}
\frac{ B^{n+1}-2B^n+B^{n-1}}{\Delta t^2} &=& -c^2\nabla_x\times (\nabla_x\times B^{n+1})\nonumber\\
\label{waveB1}
&=& -c^2\nabla_x(\nabla_x\cdot B^{n+1}) +c^2\Delta B^{n+1}. 
\end{eqnarray}
Moreover, the divergence constraint for $B$ is propagated in time  since if $\nabla_x\cdot B^n=0$, 
then the second equation of \eqref{maxwell-linear-time} ensures $\nabla_x\cdot B^{n+1}=0$. 
As a consequence, \eqref{waveB1} reduces to the following implicit time integrator
\begin{equation}
\label{schimmel}
\frac{ B^{n+1}-2B^n+B^{n-1}}{\Delta t^2}  = c^2\Delta B^{n+1}. 
\end{equation}
It is well known that this scheme is stable since the amplification factor is smaller than one. 
It means that this scheme damps high frequencies significantly compared to the exact flow. 
However, this property is essential in our case as we rely on the fact that for large $c$ 
the numerical scheme damps the high frequencies in the system to recover the correct 
asymptotic behavior. 
With the solver described in \cite{cef} (or solving \eqref{waveB1} exactly), 
this behavior is not possible.  Thus, we expect the present numerical integrator to compare unfavourably to the 
splitting described in \cite{cef} when $c$ is close to unity (even though it is a consistent numerical scheme). However, for large values of $c$ it has the decisive advantage that no CFL condition is imposed for the wave equation. 

Note that several adaptations could be introduced to make the scheme symplectic if $c$ is small, by changing for instance the right-hand side of \eqref{schimmel} to $c^2( (1 - \theta) \Delta B^{n+1} + \theta \Delta B^n)$ where $\theta$ depends on $c$, and can be chosen close to $1$ for large $c$ and close to $0$ for small $c$. However, we will not consider such modifications in the present paper. 

\subsection{Time discretization of the Vlasov equation}
\label{subsec:time_integ}
Let us now focus on the kinetic part of the RVM equation. To ensure that the Poisson equation is satisfied for all times, 
the numerical scheme should satisfy the charge conservation property (see \cite{chacon, cns, valis}). 
To accomplish this we adopt a time splitting inspired from \cite{cef}: 
first, we solve the following flow
\begin{equation}
\begin{array}{l}
\displaystyle \partial_t f + \frac{p}{\gamma} \cdot \nabla_x f  = 0, \\
\displaystyle \partial_t E =  - \int_{\R^3} \frac{p}{\gamma} f(t,x,p) \dd p { + \bar{J}(t)}, \\
\partial_t B = 0.
\end{array}
\label{Hf}
\end{equation}

\noindent Second, we solve 
\begin{equation}
\begin{array}{l}
\partial_t f +E \cdot \nabla_p f  = 0, \\
\partial_t E = 0, \quad
\partial_t B = 0.
\end{array}
\label{HE}
\end{equation}

\noindent Finally, we solve 
\begin{equation}
\begin{array}{l}
\displaystyle\partial_t f + \left(\frac{p}{\gamma}\times B\right) \cdot \nabla_p f   = 0, \\
\partial_t E = 0, \quad
\partial_t B = 0.
\end{array}
\label{HB}
\end{equation}
Note that $\nabla_p \cdot \left(\frac{p}{\gamma}\times B(x)\right)=0$ so that the transport term in \eqref{HB} 
is conservative. All of these three steps can be solved exactly in time (see \cite{cef}). 

\subsection{Application to the 1+1/2 RVM and phase space integration}
\label{1p12}
Preparing for the numerical simulations conducted in section \ref{sec:numerical-results}, we detail our numerical scheme for the 1+1/2 RVM system (see also \cite{suzuki}). 
We consider the phase space $(x_1,p_1,p_2) \in X\times \R^2$, where $X$ is a one-dimensional torus,  and the unknown functions are 
$f(t,x_1,p_1,p_2)$, $B(t,x_1)$ and $E(t,x_1) = (E_1,E_2)(t,x_1)$ which are determined by solving the following system of evolution equations
\begin{equation}
\label{eq:vmreduced}
\begin{array}{l}
\displaystyle \partial_t f + \frac{p_1}{\gamma}  \partial_{x_1} f + E \cdot \nabla_p f  + \frac{1}{\gamma} B \mathcal{J} p \cdot \nabla_p f = 0,\\
\partial_t B = -  \partial_{x_1} E_2; \\
\partial_t E_2 = -  c^2\partial_{x_1} B - \displaystyle \int_{\R^2} \frac{p_2}{\gamma} f(t,x_1,p) \dd p  + \bar{J}_2(t),\\
\partial_t E_1 = - \displaystyle \int_{\R^2} \frac{p_1}{\gamma} f(t,x_1, p) \dd p  + \bar{J}_1(t),\\
\end{array}
\end{equation}
where $p=(p_1, p_2)$, $\gamma=\sqrt{1+(p_1^2+p_2^2)/c^2}$ and  
$$
\bar{J}_i(t)=\frac{1}{|X|} \int_{X} \int_{\R^2} \frac{p_i}{\gamma} f(t,x_1,p) \dd x_1 \dd p,\quad  i=1, 2, 
$$ 
with $|X|$ the total measure of $X$; finally, $\mathcal{J}$ denotes the symplectic matrix 
$$
\mathcal{J} = \begin{pmatrix}
0 & 1 \\Ê-1 & 0
\end{pmatrix}.
$$
This reduced system corresponds to choosing an initial value of the form 
$$
E(x_1,x_2,x_3) = \begin{pmatrix} E_1(x_1) \\ E_2(x_1) \\ 0 \end{pmatrix}
\quad \mbox{and}\quad B(x_1,x_2,x_3) = \begin{pmatrix} 0 \\ 0 \\ÊB(x_1) \end{pmatrix}
$$
and a $f$ depending on $x_1$ and $(p_1,p_2)$ only in the system \eqref{eq:vm}. Then it can be easily checked that this structure is preserved by the exact flow. We refer the reader to \cite{chengy, suzuki} for more details.

Let us now consider the splitting scheme introduced in the previous subsections in more detail in the context of the 1+1/2 RVM model.
We denote by $f^n$, $E_1^n$, $E_2^n$ and $B^n$ approximations of the exact solution at time $t_n = n \Delta t$.

\subsubsection{First step}
The first step of the splitting consists in advancing \eqref{Hf} in time which, in our 1+1/2 RVM framework, can be written in Fourier space as follows 
\begin{equation}
\label{Hf_reduced}
\partial_t \hat{f} + \frac{p_1}{\gamma} i k \hat{f} = 0, \;\; \partial_t \hat{E}_1 = -\int_{\R^2} \frac{p_1}{\gamma} \hat{f} \dd p + {\bar{J_1}},\;\; \partial_t \hat{E}_2 = -\int_{\R^2} \frac{p_2}{\gamma} \hat{f} \dd p + {\bar{J_2}},
\end{equation}
where $\,\hat{}\,$ denotes the Fourier transform in the spatial variable only. We extend \cite{cef} to
the relativistic case: first, $\hat{f}$ can be computed exactly from $\hat{f}^n$ by integrating directly between $0$ and $\Delta t$
$$
\forall\, k \in \Z, \quad 
\hat{f}^{\star} = \hat{f}^n \exp(- i p_1  k \Delta t/\gamma). 
$$
Then, the equation for $E_1$ can be solved exactly in time (for $k\neq 0$)
\begin{eqnarray*}
\hat{E}_1^{\star} &=& \hat{E}_1^n - \int_{\R^2} \frac{p_1}{\gamma} \int_0^{\Delta t} \hat{f}(t)  \dd t  \dd p \nonumber\\
&=& \hat{E}_1^n - \int_{\R^2} \frac{p_1}{\gamma} \hat{f}^n  \int_0^{\Delta t} \exp(- i p_1  k (t-t^n)/\gamma) \dd t  \dd p \nonumber\\
&=& \hat{E}_1^n - \int_{\R^2} \frac{p_1}{\gamma} \hat{f}^n \left[   \frac{-1}{i k p_1/\gamma} (\exp(-ik p_1 \Delta t/\gamma)-1 )  \right] \dd p \nonumber\\
&=&  \hat{E}_1^n + \frac{1}{ik} \left[ \int_{\R^2} (\hat{f}^\star - \hat{f}^n )\dd p\right].
\end{eqnarray*}
The same procedure can be applied to the equation for $E_2$ 
\begin{eqnarray*}
\hat{E}_2^{\star} &=& \hat{E}_2^n - \int_{\R^2} \frac{p_2}{\gamma} \hat{f}^n \left[   \frac{-1}{i k p_1/\gamma} (\exp(-ik p_1 \Delta t/\gamma)-1 )  \right] \dd p \nonumber\\
 &=& \hat{E}_2^n +  \frac{1}{i k} \int_{\R^2} \frac{p_2}{p_1} \hat{f}^n \left[\exp(-ik p_1 \Delta t/\gamma)-1   \right] \dd p. \nonumber\\
\end{eqnarray*}
Numerically, the integration with respect to $p$ is done by standard quadrature formulas.

\subsubsection{Second step}
In the second step we approximate the linear part of Maxwell's equations  \eqref{maxwell-linear-time}. 
For the 1+1/2 RVM case we get (in Fourier space) 
\begin{equation}
\label{maxwell_reduced}
\partial_t \hat{B} = -  ik \hat{E}_2, \;\; 
\partial_t \hat{E}_2 = -  c^2 ik \hat{B}, \;\; 
\partial_t \hat{E}_1 = 0, 
\end{equation}
with the initial conditions $\hat{B}=\hat{B}^n$, $\hat{E}_1=\hat{E}_1^\star$, $\hat{E}_2=\hat{E}_2^\star$  
($\hat{E}_1^\star$ and $\hat{E}_2^\star$ are computed in the last step). 
The use of  an implicit Euler scheme in time to ensure stability with respect to $c$ yields the formula
\begin{eqnarray*}
\hat{E}_2^{n+1} &=& \hat{E}_2^{n} - c^2 \Delta t\,  i k \hat{B}^{n+1} \nonumber\\
\hat{B}^{n+1} &=& \hat{B}^{n} - \Delta t\,  i k \hat{E}_2^{n+1}.  %
\end{eqnarray*}
Note that $\hat{E}_1$ is unchanged and thus $\hat{E}_1^{n+1} = \hat{E}_1^{\star}$. 
These equations can be cast into the following $2$x$2$ matrix system 
\begin{equation}
\label{linear_maxwell_euler}
\left(
\begin{array}{ll}
\hat{E}_2^{n+1} \\
\hat{B}^{n+1}
\end{array}
\right)
=
\frac{1}{1+\Delta t^2 c^2 k^2}
\left(
\begin{array}{ll}
1 & -c^2 \Delta tik \\
-\Delta t ik & 1
\end{array}
\right)
\left(
\begin{array}{ll}
\hat{E}_2^{\star} \\
\hat{B}^{n}
\end{array}
\right). 
\end{equation}

\subsubsection{Third step}
In the third step we solve \eqref{HE}. Using the electric field $E^{n+1}$ computed in the previous step, it becomes 
\begin{equation}
\label{tp1}
\partial_t f + E^{n+1}\cdot \nabla_p f = 0. 
\end{equation}
As the electric field is kept constant during this step, the solution of this equation is explicitly given by 
$$
f^{\star\star}(x,p) = f^{\star}(x, p-\Delta t E^{n+1}).  
$$ 
The evaluation of $f^\star$ at the point $p-\Delta t E^{n+1}$ is performed using a 2-dimensional 
interpolation (using Lagrange interpolation of degree $3$). 

\subsubsection{Fourth step}
In this last step, we solve \eqref{HB}. Using the magnetic field $B^{n+1}$ computed in the second step, we have to solve 
\begin{equation}
\label{tp2}
\partial_t f + \frac{1}{\gamma} B^{n+1}{\mathcal J} p\cdot \nabla_p f = 0. 
\end{equation}
The solution of this equation can be written as follows
\begin{equation}
\label{interp_2d_p}
f^{n+1}(x,p) = f^{\star\star}(x, P(t^n; t^{n+1}, p)), 
\end{equation}
where $P(t^n; t^{n+1}, p)$ is the solution at time $t^n$ of the characteristics equation taking the value $p$ at time $t^{n+1}$, i.e.
$$
\frac{d P}{dt} = \frac{1}{\sqrt{1+|P|^2/c^2}}B^{n+1} {\mathcal J} P, \;\;\;\; P(t^{n+1}) = p, \;\;\;\; t\in [t^n, t^{n+1}].  
$$
This ordinary differential equation can be solved analytically since $\gamma$ is constant on each trajectory and $B^{n+1}$ is independent of $p$. Then, a 2-dimensional interpolation (using Lagrange interpolation of degree $3$) is performed in \eqref{interp_2d_p} in order to compute $f^{n+1}$.  In the non-relativistic case, this step is simpler since a directional splitting reduces the problem to a sequence of one-dimensional transport equations.

\subsubsection{Algorithm}
We summarize the main point of the proposed algorithm, starting from $f^n, E_1^n, E_2^n, B^n$: 
\begin{itemize}
\item compute $f^\star, E_1^\star, E_2^\star$ from $f^n, E_1^n, E_2^n$ by solving \eqref{Hf_reduced} with step size $\Delta t$, 
\item compute $E_1^{n+1}, E_2^{n+1}, B^{n+1}$ from $E_1^\star, E_2^\star, B^n$ by solving \eqref{maxwell_reduced} using the implicit Euler method with step size $\Delta t$, 
\item compute $f^{\star \star}$ from $f^\star$ by solving \eqref{tp1} with step size $\Delta t$, %
\item compute $f^{n+1}$ from $f^{\star \star}$ by solving \eqref{tp2} with step size $\Delta t$. %
\end{itemize}

\subsubsection{Asymptotic preserving property}
We are interested here in the asymptotic behavior of the proposed numerical 
scheme when $c$ goes to $+\infty$, for a fixed time step $\Delta t$ and  independently from the initial condition. 

From the second step, we immediately get from \eqref{linear_maxwell_euler} 
that the magnetic field $\hat{B}^{n+1}$ goes to zero when $c\to +\infty$, for all $n\geq 0$. 
Moreover, again from \eqref{linear_maxwell_euler}, 
the term $\hat{E}_2^{n+1}$ goes to $-i/(k\Delta t) \hat{B}^n$ when  $c\to +\infty$. 
Hence we get that $\hat{E}_2^{n+1}$ goes to zero 
when $c\to +\infty$, for all $n\geq 1$. Note that even if the initial condition is not consistent with the asymptotic behavior 
(i.e. $\hat{E}_2^{0}\neq 0$ or $\hat{B}^0\neq 0$), the numerical scheme we propose imposes, after the first steps, that 
$\hat{E}_2^{2}$ and $\hat{B}^1$ become small as $c\to +\infty$. This is related to the strong asymptotic property 
which does not require that the initial data are well-prepared, typically $\hat{B}^0=\mathcal{O}(1/c)$.  
Thus, the only field that does not vanish when $c\to +\infty$ is the electric field $\hat{E}_1$.   

Since we have ensured that $E_2$ goes to zero as $c$ goes to $+\infty$, the third step reduces 
to a one-dimensional transport in the $p_1$ direction due to the effect of $E_1^{n+1}$ 
(which has been computed in the first step). 
Similarly, since $B$ goes to zero as $c$ goes to $+\infty$, the last step 
leaves $f$ unchanged. 

The numerical method described is first order in time for a fixed value of $c$, 
and satisfies the asymptotic preserving property. More precisely, as $c$ goes to infinity and for a fixed $\Delta t$, the algorithm reduces to
\begin{itemize}
\item solve $\partial_t \hat{f} + p_1 i k \hat{f} = 0$ and $\partial_t \hat{E}_1 = -\int_{\R^2} p_1 \hat{f} \dd p + \bar{J_1}$ with step size $\Delta t$. This gives
$$
\hat{f}^{\star}(k,p) = Ê\hat{f}^n  \exp(- i p_1  k \Delta t),
$$
and
$$
 \hat{E}_1^{n+1}(k) = \hat{E}_1^n(k)+\frac{1}{ik} \left[ \int_{\R^2} (\hat{f}^\star(k,p) - \hat{f}^n(k,p) )\dd p\right].
$$
\item solve $\partial_t f + E_1\partial_{p_1} f =0$ with step size $\Delta t$ and the $E_1$ computed in the previous step. This gives
$$
f^{n+1}(x,p) = f^{\star}(x, p_1-\Delta t E_1^{n+1}, p_2).  
$$
\end{itemize}
This so-obtained asymptotic numerical scheme corresponds to a two-term splitting 
for the one-dimensional Vlasov--Amp\`ere equations in the variables $(x_1, p_1)$. 
We emphasize that this scheme is consistent with the continuous asymptotic one-dimensional 
Vlasov--Poisson model since it preserves the charge exactly (see \cite{cef}); indeed, if the Poisson 
equation is satisfied initially, it is satisfied for all time due to the fact that we solve Amp\`ere's equation exactly.

\subsection{Extension to second order\label{sec:2nd}}

{\color{black}
In practical simulations constructing a scheme that is at least of second order is a necessity in order to obtain good accuracy. 
The previous scheme can be easily extended to second order by using the symmetric Strang splitting. In addition, higher order splitting methods can easily be constructed by composition (see \cite{hlw}). The missing crucial ingredient, however, is an integrator for the linear part of Maxwell's equations \eqref{linear_maxwell} that is of the appropriate order (so far we have only considered the first order implicit Euler scheme).

If we consider the problem of solving the relativistic Vlasov--Maxwell system on a tensor product domain, it is even possible to exactly integrate the linear part of Maxwell's equations. This is possible since \eqref{linear_maxwell} is a linear system with constant coefficients and therefore all the Fourier modes decouple (see equation \eqref{maxwell_reduced}). The resulting 2x2 (complex) matrix exponential can be computed analytically. However, as is evident from Figure \ref{fig:ldmag-stability} (top-left), using this approach we do not even converge towards the correct limit (i.e.~we do not observe the correct Landau damping rate for large values of $c$). This is due to the fact that the exact solution does not damp high frequencies at all. Consequently, we do not approach the electrostatic limit which the RVM system only attains in a weak sense.

\begin{figure}
	\centering
	\includegraphics[width=12cm]{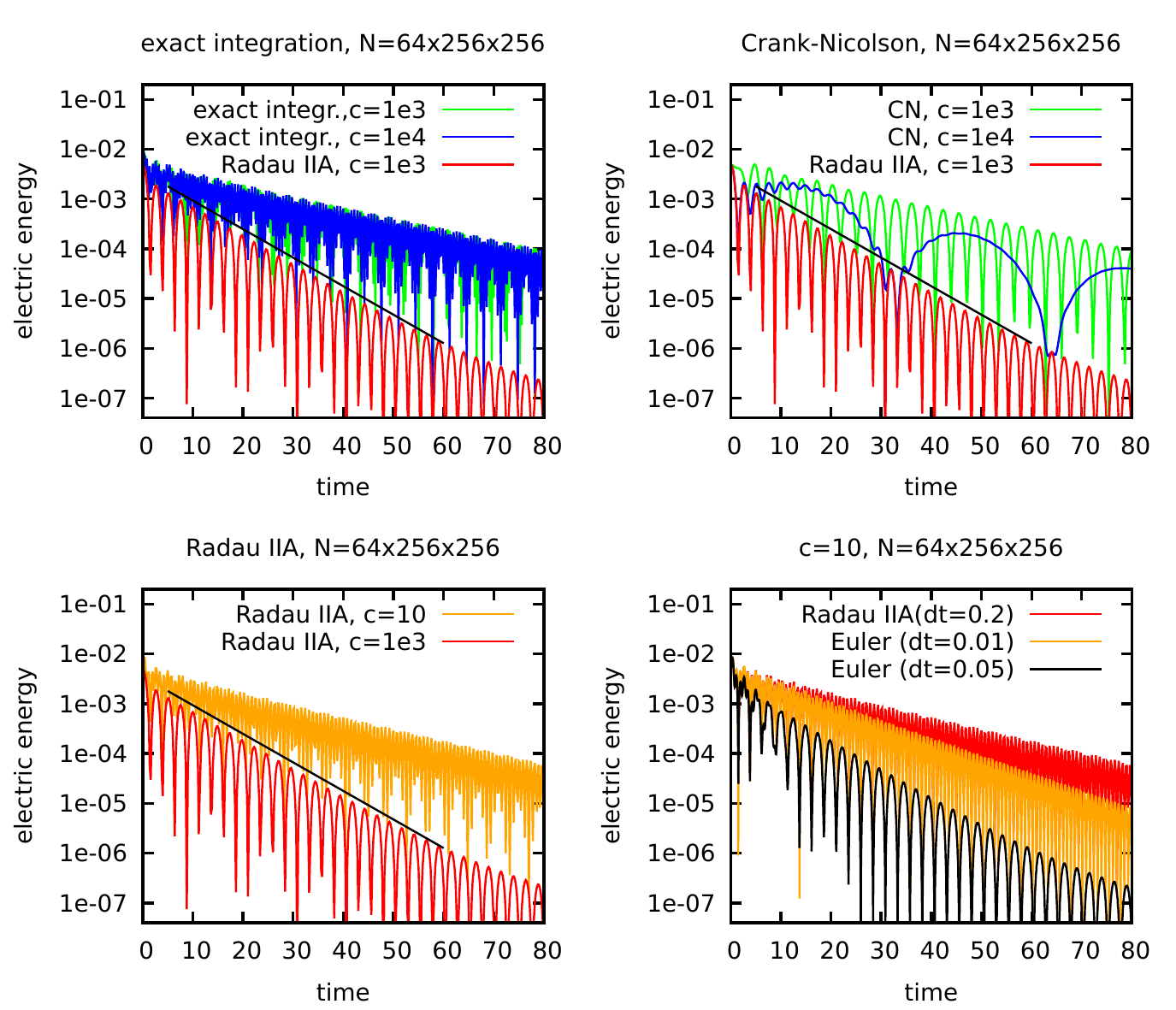}
	\caption{\color{black}The Landau damping problem for a number of different numerical schemes is shown. In the top row exact integration and the Crank--Nicolson method are used. In this case we do not observe convergence to the correct limit. At the bottom-left graph we see that using the third order Radau IIA method results in the correct limit behavior. The bottom-right graph compares the first order implicit Euler method with the third order Radau IIA method for $c=10$. The analytically derived decay rate for the Vlasov--Poisson system is shown as a black line. In all simulations, if not indicated otherwise, the time step size is chosen as $\Delta t = 0.2$. \label{fig:ldmag-stability}}
\end{figure}

The Crank--Nicolson method is a very widely used numerical scheme that is second order accurate and only requires the inversion of a single system of linear equations. It thus incurs the same computational cost as the implicit Euler method. Unfortunately, it suffers from the same shortcoming as integrating the linear part of Maxwell's equations exactly (i.e.~we do not observe the correct limit for large values of $c$). The corresponding numerical results are shown in Figure \ref{fig:ldmag-stability} (top-right).

We can perform a linear analysis of these methods by recognizing that equation (\ref{maxwell_reduced}) can be diagonalized. The corresponding eigenvalues are $\pm i c k$. Therefore, it is sufficient to only consider the stability function $\phi(z)$ of the numerical integrators used (these are listed in Table \ref{tab:stability-functions}). In particular, 
the exact integration and the Crank--Nicolson method give $\lim_{z\to\infty} \phi(z) = 1$ which means that 
there is no damping for $c \to +\infty$. Therefore using these methods our numerical scheme does not converge to the classical limit.

\begin{table}
	\centering
\textcolor{black}{
	\begin{tabular}{c|ccccc}
		Method & exact & Crank--Nicolson & imp. Euler & Radau & SDIRK \\
		\hline
		$\phi(z)$ & $\mathrm{e}^{z}$ & $\frac{1+\frac{z}{2}}{1-\frac{z}{2}}$ & $\frac{1}{1-z}$ & $\frac{1+ \frac{1}{3}z}{1-\frac{2}{3}z+\tfrac{1}{6}z^{2}}$ & $\frac{1+(\sqrt{2}-1)z}{(1+(\frac{\sqrt{2}}{2}-1)z)^{2}}$ \\
	\end{tabular}
	\caption{The stability function $\phi(z)$ for the exact integration, the Crank--Nicolson method, the implicit Euler method, the Radau IIA method of order three, and the L-stable SDIRK method given in Figure \ref{fig:ldmag-lstable} are listed. \label{tab:stability-functions}}
}
\end{table}

Implicit Runge--Kutta methods have been constructed so as to satisfy $\lim_{z\to\infty} \phi(z)=0$. These so-called $L$-stable methods\footnote{A method is called $L$-stable if it is $A$-stable and satisfies $\lim_{z\to\infty} \phi(z)=0$.} have the property that for $c\to +\infty$ the magnetic field vanishes after a single time step. This is clearly a desirable property in the present situation. The most commonly used member of this class are the so-called Radau IIA methods. The Radau IIA method with $s$ stages converges with order $2s-1$. In fact, the implicit Euler method is identical to the Radau IIA method with $s=1$.

The numerical results for the Strang splitting using the Radau IIA method for the linear part of Maxwell's equations are 
 shown in Figure \ref{fig:ldmag-stability} (bottom-left). In this case we do observe the correct behavior in the classical limit (as is demonstrated by the comparison to an analytically derived result for the Vlasov--Poisson system). Let us also note that, while the implicit Euler method is able to recover the correct limit, it creates a numerical damping even for relatively small values of $c$. This is a numerical artefact that vanishes as we decrease the time step size. However, in order to obtain results comparable to the third order Radau IIA method the implicit Euler method has to use a time step size that is at least a factor of 20 times smaller (see Figure \ref{fig:ldmag-stability} bottom-right; in both cases we use the second order Strang splitting scheme for the whole RVM system). Therefore, we will use the third order Radau IIA method in all the simulations that have been conducted in this paper.

A disadvantage of the Radau IIA family of methods is that they are fully implicit (see the Butcher Tableau in Table \ref{coef_radau}). In general, we thus have to solve a nonlinear system of equations coupling all stages of the numerical method. In the present case this is not a severe restriction for the following reasons.  First, once we apply the splitting scheme, the resulting Maxwell's equations are linear 
(see  equation \eqref{linear_maxwell}). Thus, we apply the Radau IIA method to a linear system and no Newton iteration is required. Second, and most important, in Fourier space the different modes decouple. Thus, for a Fourier based space discretization the Radau IIA method of third order yields a complex $4$x$4$ system of linear equations for each mode. This system can be solved analytically. The resulting expression is employed in our implementation. However, for applications where Fourier techniques are not applicable, inverting the linear system required to evaluate these implicit methods can incur a significant computational cost. In this case we can either use the techniques described in \cite[Chapt. IV.8]{hwstiff} to reduce the size of the linear system for the Radau IIA family of methods or employ a single diagonally implicit (SDIRK) $L$-stable method. Second order $L$-stable SDIRK methods with two stages have been constructed and we have implemented one such method. The corresponding numerical results, which show convergence to the correct limit as $c\to +\infty$, are displayed in Figure \ref{fig:ldmag-lstable}. 

Even though this $L$-stable SDIRK method is a viable alternative to the Radau IIA method discussed earlier, in the present paper we will only report results using the latter scheme. This is due to the fact that the Radau IIA method is accurate to third order and there is no additional computational cost as we exclusively use Fourier techniques in order to discretize space.

\begin{figure}\color{black}
	\begin{center}
	\begin{minipage}{0.4\textwidth}
		\centering
	\begin{tabular}{c|cc}
	$\gamma$ & $\gamma$ & $0$\\
	&&\\
	$1-\gamma$ & $1-2\gamma$ &  $\gamma$\\
	\hline\\
	&$1/2$ & $1/2$
	\end{tabular}
	\end{minipage}
	\begin{minipage}{0.4\textwidth}
		\centering
	\includegraphics[width=6cm]{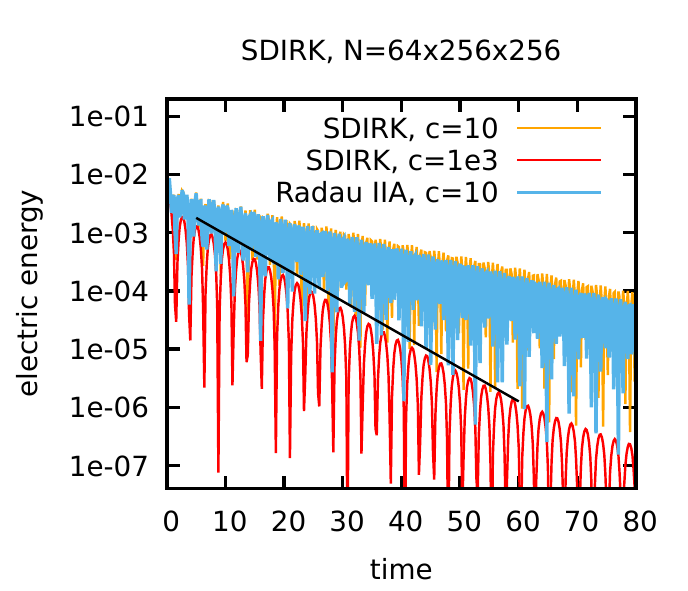}
	\end{minipage}
	\end{center}

	\caption{\color{black}The Butcher Tableau for the second order SDIRK method with two stages is shown on the left ($\gamma=1-1/\sqrt{2}$). The figure on the right shows that for large $c$ we recover the classical limit when using this method to solve the linear part of Maxwell's equations. The analytically derived decay rate for the Vlasov--Poisson system is shown as a black line. In all simulations the time step size is chosen as $\Delta t = 0.2$.\label{fig:ldmag-lstable}}
\end{figure}

Now let us proceed by describing the second order numerical scheme that is used in all simulations in the next section:
}

\begin{itemize}
	\item compute $f^\star, E^\star$ from $f^n, E^n$ by solving \eqref{Hf} with step size $\tfrac{1}{2} \Delta t$, 
	\item compute $E^{n+1/2}, B^{n+1/2}$ from $E^\star, B^n$ by solving \eqref{linear_maxwell_euler} %
	using the third order Radau IIA method with step size $\tfrac{1}{2}\Delta t$, 
	\item compute $f^{\star \star}$ from $f^\star$ by solving \eqref{HE} with $E=E^{n+1/2}$ with step size $\frac{1}{2}\Delta t$, 
	\item compute $f^{\star \star \star}$ from $f^{\star \star}$ by solving \eqref{HB} with $B=B^{n+1/2}$ with step size $\Delta t$,
	\item compute $f^{\star \star \star\star}$ from $f^{\star \star \star}$ by solving \eqref{HE} with $E=E^{n+1/2}$ with step size $\frac{1}{2}\Delta t$, 
	\item compute $E^{\star \star}, B^{n+1}$ from $E^{n+1/2}, B^{n+1/2}$ by solving \eqref{linear_maxwell_euler} %
	using the third order Radau IIA method with step size $\tfrac{1}{2}\Delta t$. 
	\item compute $f^{n+1}, E^{n+1}$ from $f^{\star \star \star\star}, E^{\star\star}$ by solving \eqref{Hf} with step size $\tfrac{1}{2} \Delta t$, 
\end{itemize}

\begin{table}
\begin{center}
\begin{tabular}{c|cc}
$1/3$ & $5/12$ & $-1/12$\\
&&\\
$1$ & $3/4$ &  $1/4$\\
\hline\\
&$3/4$ & $1/4$
\end{tabular}
\end{center}
\caption{Butcher Tableau for the Radau IIA method.}
\label{coef_radau}
\end{table}

\subsection{Extension to the general case\label{sec:dim}}
\textcolor{black}{
	In this section, we briefly discuss how to extend the algorithm described in section \ref{1p12} to the general case of three spatial and three velocity directions (i.e.~the general formulation stated in equation \eqref{eq:vm}). \\
	The first step can be generalized easily since we obtain (by using a Fourier transformation with the variables $k=(k_1, k_2, k_3)$ in space)  
$$
\hat{f}^{\star} = \hat{f}^n \exp(- i p \cdot  k \Delta t/\gamma), 
\mbox{ and }  
\hat{E}^\star = \hat{E}^n + \int_{\R^3} \frac{p}{ik\cdot p} (\hat{f}^\star-\hat{f}^n)  \dd p.  
$$
Let us remark that, similar to the one-dimensional case, charge conservation is ensured in this step. \\
The second step concerns the linear part of the Maxwell equations. 
In Fourier variables, they can be written as 
$$
\partial_t \hat{E} = c^2 i k \times \hat{B}, \;\;\; \partial_t \hat{B} = - i k \times \hat{E},  
$$
and as before, by introducing  $y(t)=(\hat{E}(t), \hat{B}(t))^T \in \mathbb{R}^6$, we get  
$$
\frac{d y(t)}{dt} = i\mathcal{A} y(t), \mbox{ with } y(t)=(\hat{E}(t), \hat{B}(t))^T \in \mathbb{R}^6, 
$$
and $\mathcal{A}$ is a $6$x$6$ matrix given by 
$$
\mathcal{A} =\left( 
\begin{array}{ccccccllllll}
0 & 0 & 0 & 0 & -c^2 k_3 & c^2 k_2 \\
0 & 0 & 0 & c^2 k_3 & 0 & -c^2 k_1 \\
0 & 0 & 0 & -c^2 k_2 & c^2 k_1 & 0 \\
0 & k_3 & -k_2 & 0 & 0 & 0 \\
-k_3 & 0 & k_1 & 0 & 0 & 0 \\
k_2 & -k_1 & 0 & 0 & 0 & 0 
\end{array}
\right). 
$$
Hence, one can apply Radau IIA type methods to this linear system (note that $I+i\Delta t \mathcal{A}$ is always invertible). \\
The third and fourth steps are unchanged since in both cases, the solution can be given explicitly 
using its invariance along the characteristics.  
}

\section{Numerical results \label{sec:numerical-results}}

This section is devoted to validating the numerical scheme introduced in this paper. To do so we will present and discuss the results of a number of numerical simulations for different values of the dimensionless parameter $c$. In all the numerical simulations conducted we employ the second order scheme that is described in section \ref{sec:2nd},  which is based on Strang splitting for the Vlasov equation and the third order Radau IIA method for the linear part of Maxwell's equations. We call it AP-VM and it will be compared in the regime $c\approx 1$ with the splitting 
proposed in \cite{cef} which we call H-split. \textcolor{black}{Note that the H-split method does not conserve energy up to machine precision but the numerical results in \cite{cef} show that it is superior compared to other methods from the literature}.
In the sequel, two configurations are studied: first numerical tests are conducted in the semi-relativistic 
case, considering $\gamma=1$ in \eqref{eq:vmreduced} but with different values for $c$ in Maxwell's equations. Second, 
the fully relativistic case is tackled with $\gamma=\sqrt{1+|p|^2/c^2}$ and different values for $c$. 
For these two configurations, both Landau and Weibel type problems are considered. 

\subsection{Semi-relativistic case: $\gamma=1$}

\subsubsection{Landau type problem}
First, we consider a problem that converges to a Landau damping situation as $c$ goes to infinity: We impose the following initial value for the particle density function
\[ f_0(x,p) = \frac{1}{2 \pi} \mathrm{e}^{-\frac{1}{2}(p_1^2+p_2^2)}(1+\alpha \cos k x), \;\; x\in [0, L], \; p\in [-p_{\max}, p_{\max}]^2. \]
In the Vlasov--Poisson case we would initialize the electric field according to Gauss's law. 
However, as our goal here is to stress the classical limit regime we will initialize 
the electric and magnetic field as a plane wave where equal energy is stored in the electric and magnetic field. Thus, we impose the following initial condition
\begin{equation}
\label{init_test1} 
E_1(x) = \frac{\alpha}{k} \sin k x,
	\qquad E_2(x) = 0, 
\qquad B(x) =  \frac{\alpha}{c k} \sin k x. 
\end{equation}
It is easy to verify that Gauss's law is satisfied for the initial value. As parameters we have chosen $\alpha=0.01$, $k=0.4$, $L=2 \pi/k$, and $p_{\max}=5$.

\begin{figure}
	\centering
	\includegraphics[width=12cm]{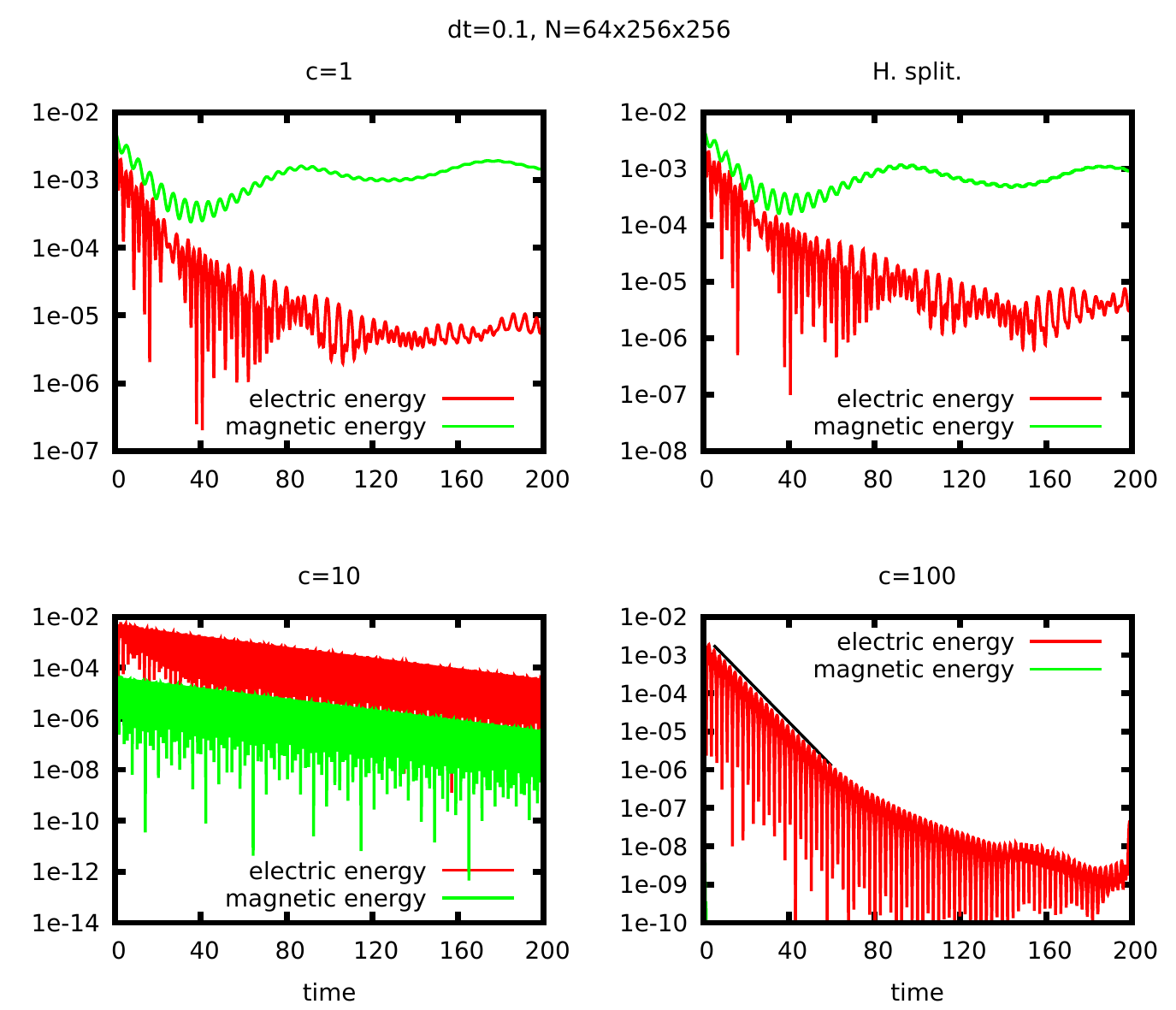}
\caption{Semi-relativistic case (Landau problem): time evolution of the electric and magnetic energy obtained by 
AP-VM (for $c=1, 10, 10^2$) and H-split ($c=1$). \label{fig:ldmag-evol}}
\end{figure}

\begin{figure}
	\centering
	\includegraphics[width=12cm]{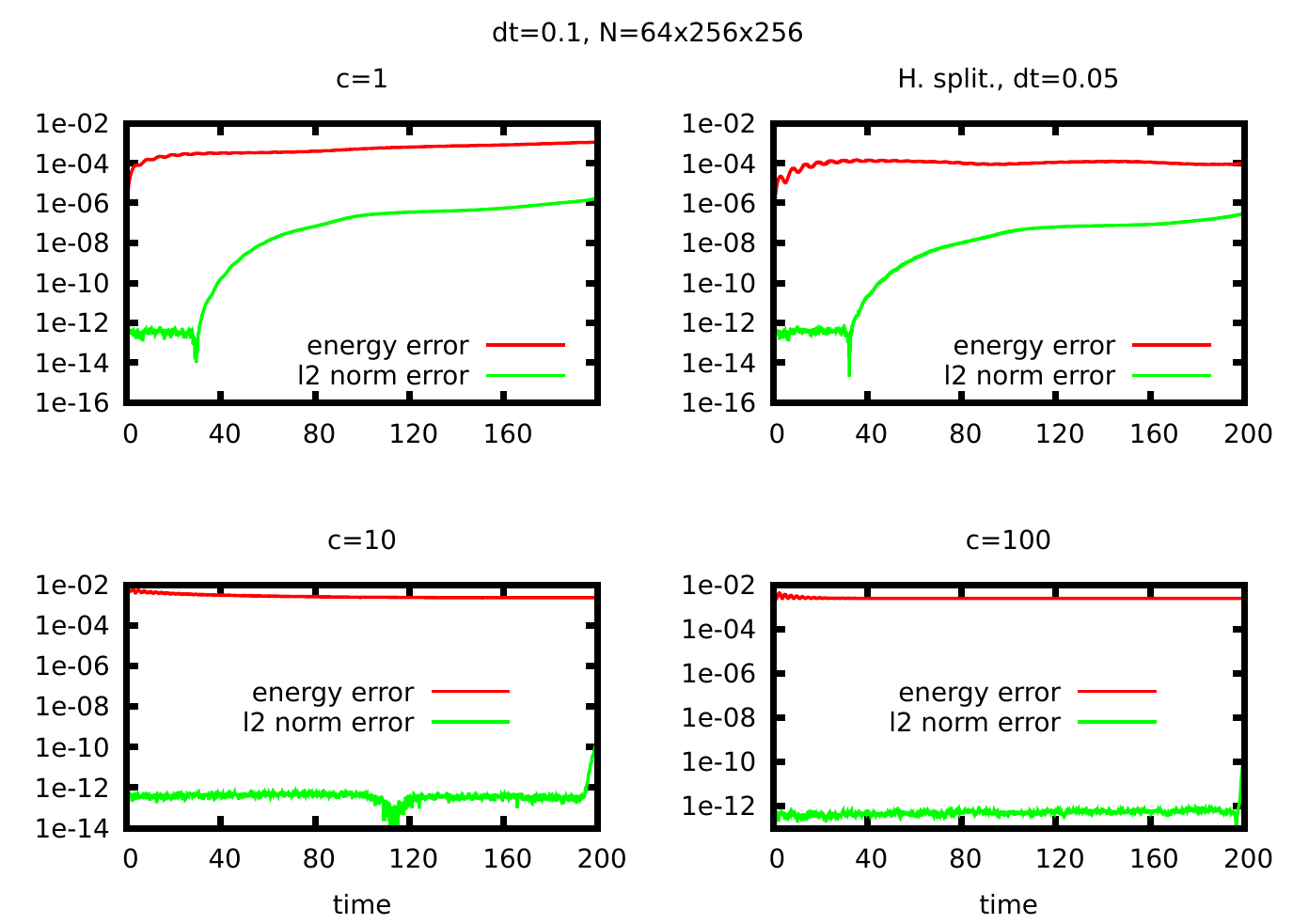}
\caption{Semi-relativistic case (Landau problem): time evolution of the error in energy and $L^2$ norm obtained by 
AP-VM  (for $c=1, 10, 10^2$) and H-split ($c=1$).\label{eq:ldmag-energ}}
\end{figure}

The numerical results are shown in Figure \ref{fig:ldmag-evol} where the time evolution of the 
electric and magnetic energies (given by $\Hc_E$ and $\Hc_B$ in \eqref{ham-split}) are 
shown for different values of $c$ ($c=1, 10, 10^2$), 
with a fixed set of numerical parameters $\Delta t=0.1, N_{x_1}=64, N_{p_1}=N_{p_2}=256$. 
We also plot the results obtained by H-split (proposed in \cite{cef}) for $c=1$ and $\Delta t=0.05$ 
in order to compare with AP-VM. It appears that AP-VM behaves very well in this regime. 
Moreover, when $c$ is large, we observe excellent agreement with the analytic results for Landau damping 
rate (the theoretical damping rate of the black line is $-0.0661$). 
The same comments apply with respect to the time evolution 
of the error in energy (defined as $|\Hc(t)-\Hc(0)|$, where $\Hc$ is defined by \eqref{ham-split}) 
and the error in the relative $L^2$ norm (in $x$ and $p$) of $f$ (defined as $\|f(t)-f(0)\|_{L^2}$):  these two quantities (which 
are preserved in time) are shown in Figure \ref{eq:ldmag-energ}.  
Indeed, when $c=1$, H-split (with $\Delta t=0.05$) and AP-VM show similar behavior.
Moreover, we can observe that the $L^2$ norm of $f$ is very well preserved when $c$ becomes large. 
This might be due to the use of Fourier methods to approximate the transport operators. 

We also look at the error in $L^\infty$ norm of the difference between the different unknown of the Vlasov--Maxwell 
system at a given $c$ ($f^c, E_1^c, E_2^c, B^c$) and the unknown of the asymptotic Vlasov--Poisson model 
($f^\infty, E_1^\infty, E_2^\infty=0, B^\infty=0$). It is known from \cite{schaeffer_cmp86} that 
(in the fully relativistic case), this error is bounded by $c^{-1}$ (with well-prepared initial data). The results we obtained 
with the initial data \eqref{init_test1} are given 
in Table \ref{table_error}. It appears that for $E_1$, the rate is stronger (the machine precision is fastly 
reached so that the last two rates are not very meaningful), for $E_2$, the rate is about $2$, for $B$ the rate 
is about $2$ 
(which corresponds to a rate of $1$ in the scaling used in \cite{schaeffer_cmp86}), 
and for $f$, the rate is about $2$.%

\begin{table}
\begin{tabular}{|c||c|c||c|c||c|c||c|c|}
\hline
$c$ &      $E_1$ error  &  rate  &        $E_2$ error   & rate &         $B$ error   & rate &         $f$ error  & rate \\
 \hline\hline
1 &   2.55e-04     &   -  &    3e-03     &   - & 6.05e-03     &   -  &    4.25e-04    &    -   \\
    \hline
 5 &        2.87e-06 &      -2.79 &        1.47e-02 &       1.15 &      2.49e-04 &      -1.98 &        1.56e-05 &      -2.05 \\
     \hline
25 &   3.78e-08  &-5.48 &    2.22e-03 &-0.02 & 2.97e-04 &-1.87 &    1.42e-06 &-3.54 \\
    \hline
125 &   5.11e-12 &-5.53 &    5.81e-07 &-5.12 &2.42e-09 &-7.28 &    2.15e-07 &-1.17 \\
    \hline
625 &   1.53e-14 &-3.61 &    2.58e-08 &-1.93 & 5.07e-11 &-2.40 &    9.8e-09 &-1.92 \\
    \hline
3125 &   1.19e-15 &-1.59 &    2.8e-10 &-2.81 &6.19e-13 &-2.74 &    4.87e-10 &-1.86 \\
    \hline
\end{tabular}
\medskip
\medskip
\caption{This table shows the difference between the numerical solution of the Vlasov--Maxwell system and the asymptotic Vlasov-Poisson system as a function of $c$. }
\label{table_error}
\end{table}

\subsubsection{Weibel type problem}
Next we consider the so-called Weibel instability. The Weibel instability is present in plasma systems with a temperature anisotropy. A small perturbation in such a system leads to an exponential growth in the magnitude of the magnetic field. The growth in amplitude eventually saturates due to nonlinear effects. The Weibel instability is considered a challenging problem for numerical simulations and is therefore often used as a test case for Vlasov--Maxwell solvers (see \cite{cpbm, cglm, cef, pcp, suzuki}).
Here we impose the following initial conditions for the particle density
\begin{equation} 
\label{init_weibel}
f_0(x,p)=\frac{1}{\pi p_{th}^2 \sqrt{T_r}} \mathrm{e}^{-(p_1^2+p_2^2/T_r)/p_{th}^2}(1+\alpha \cos k x), \; x\in [0, L], p\in [-p_{\max}, p_{\max}]^2,  
\end{equation} 
and the field variables
\[ E_1(x)=\frac{\alpha}{k} \sin kx, \qquad
	E_2(x)=0, \qquad
B(x)=\frac{\alpha}{c k} \cos k x. \]
As parameters we have chosen $\alpha=10^{-4}$, $k=1.25$, $T_r=12$, $p_{th}=0.02$, $L=2 \pi / k$, and $p_{\max}=0.3$. 
We compare the results obtained by AP-VM and by H-split.  %

\begin{figure}
	\centering
	\includegraphics[width=12cm]{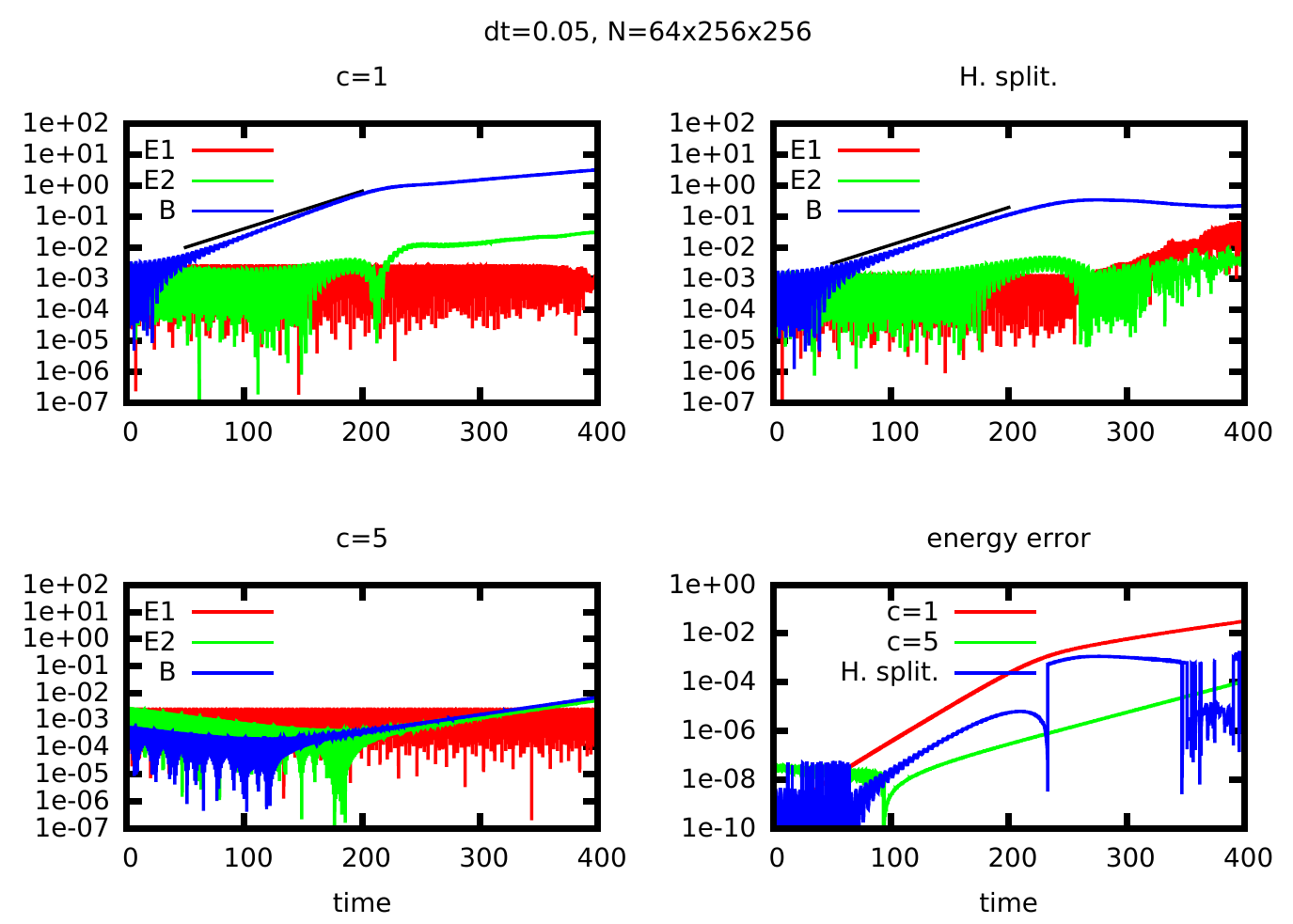}
\caption{Semi-relativistic case (Weibel problem): time evolution of the most unstable mode ($k=1.25$) 
of the magnetic and the two electric fields. Top left: $c=1$, AP-VM. Top right: $c=1$, H-split. 
Bottom left: $c=5$, AP-VM. Bottom right: time evolution of the energy error for AP-VM ($c=1, 5$) 
and for H-split ($c=1$). \label{eq:tt-evol}}
\end{figure}

\begin{figure}
	\centering
	\includegraphics[width=12cm]{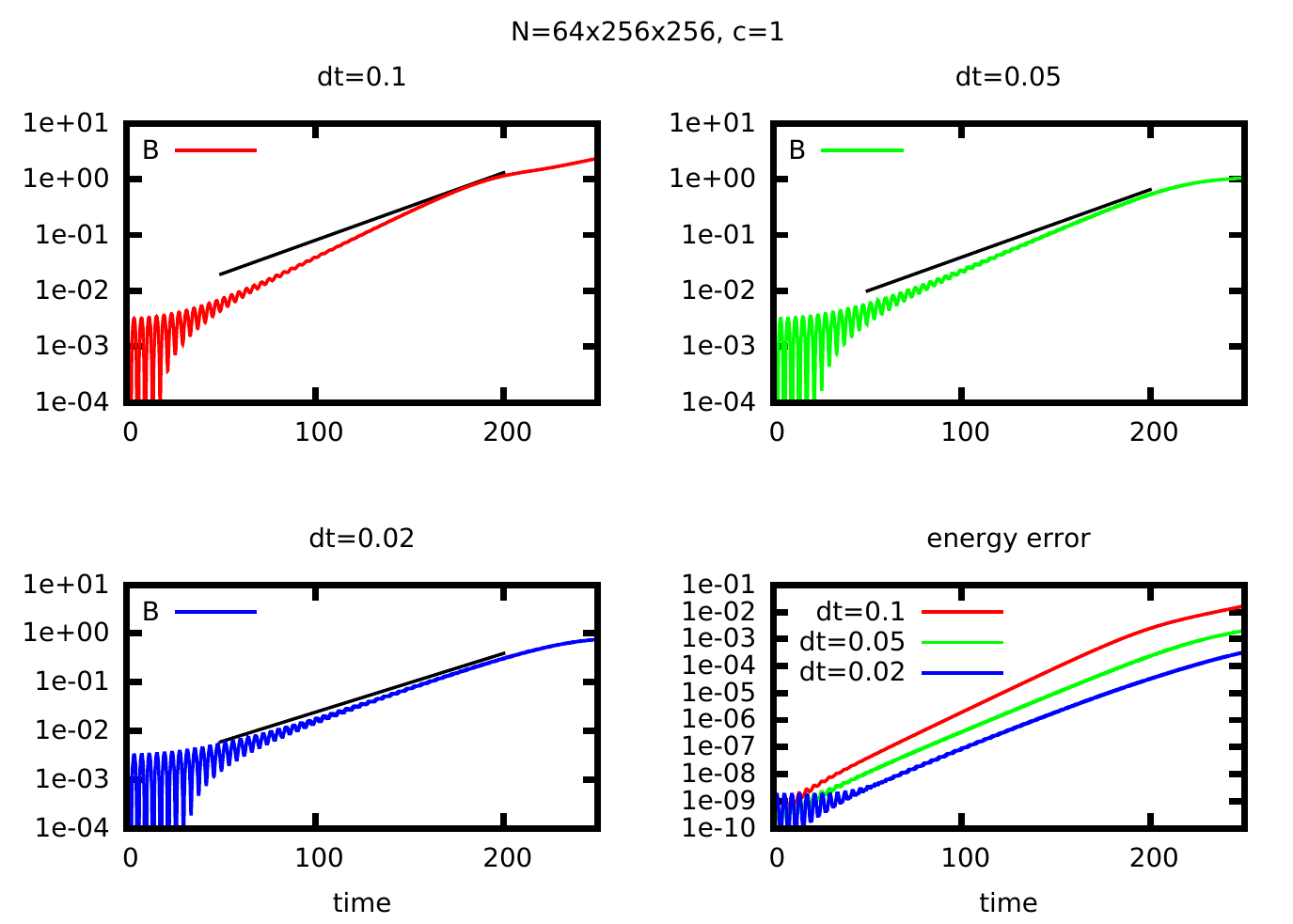}
\caption{Semi-relativistic case (Weibel problem): time evolution of the most unstable mode of the magnetic field 
(and the corresponding theoretical growth rate) obtained by the asymptotic preserving scheme 
for different time steps $\Delta t=0.1, 0.05, 0.02$. Bottom right: time evolution of the energy error for different time steps. \label{eq:tt-growth}}
\end{figure}

We are interested in the time evolution of the most unstable Fourier mode (namely $k=1.25$) of the electric 
and magnetic fields $E_1, E_2, B$, and in the time evolution of the relative total energy $\Hc(t)-\Hc(0)$. 
The numerical results are shown in Figures \ref{eq:tt-evol} and \ref{eq:tt-growth}. 
We observe that for the time step chosen $\Delta t=0.05$, H-split gives significantly better agreement 
with the growth rate derived in subsection \ref{appendix3} compared to AP-VM. Note however that it is entirely expected that preserving the Hamiltonian structure gives better qualitative agreement with the exact solution. In addition, the diffusive nature of the Radau method employed introduces significant errors in the case where $c=1$ (see the time evolution of the total 
energy). Despite this, the linear phase is well reproduced. Moreover, let us note that AP-VM is consistent since it converges 
(when $\Delta t$ is decreased sufficiently) to the correct behavior, as can be observed from Figure \ref{eq:tt-growth}. 
This also enables us to check that our AP-VM scheme is second order in time. 

As we increase the dimensionless parameter $c$ we expect the Weibel instability to cease. On physical grounds one would argue that the instability cannot exist in the electrostatic regime as the Vlasov--Poisson system does not include any magnetic effects. This is confirmed by the linear analysis that has been conducted 
(in section \ref{sec:dispersion}) which shows that even for moderate values of $c$ no unstable magnetic modes exist. The test, for AP-VM, is then to work well in this limit. We observe from Figure \ref{eq:tt-evol} that 
the energy conservation improves dramatically as $c$ increases.
For any value of $c$ larger than $5$ no instability can be observed in the case of the asymptotic scheme, and the total energy is well preserved. 
Let us note that due to the CFL restriction for the integration of the field variables, the scheme H-split is forced to take excessively small step sizes as $c$ increases.  For small $c$ this can be alleviated to some extend by performing substepping for Maxwell's equations (as pointed out in \cite{cef}); however, for medium to large $c$, H-split is computationally infeasible. On the other hand, AP-VM is unconditionally stable so that it does not suffer from this step size restriction. 

\subsection{Fully relativistic case}

\subsubsection{Landau type problem}
\begin{figure}
	\centering
	\includegraphics[width=12cm]{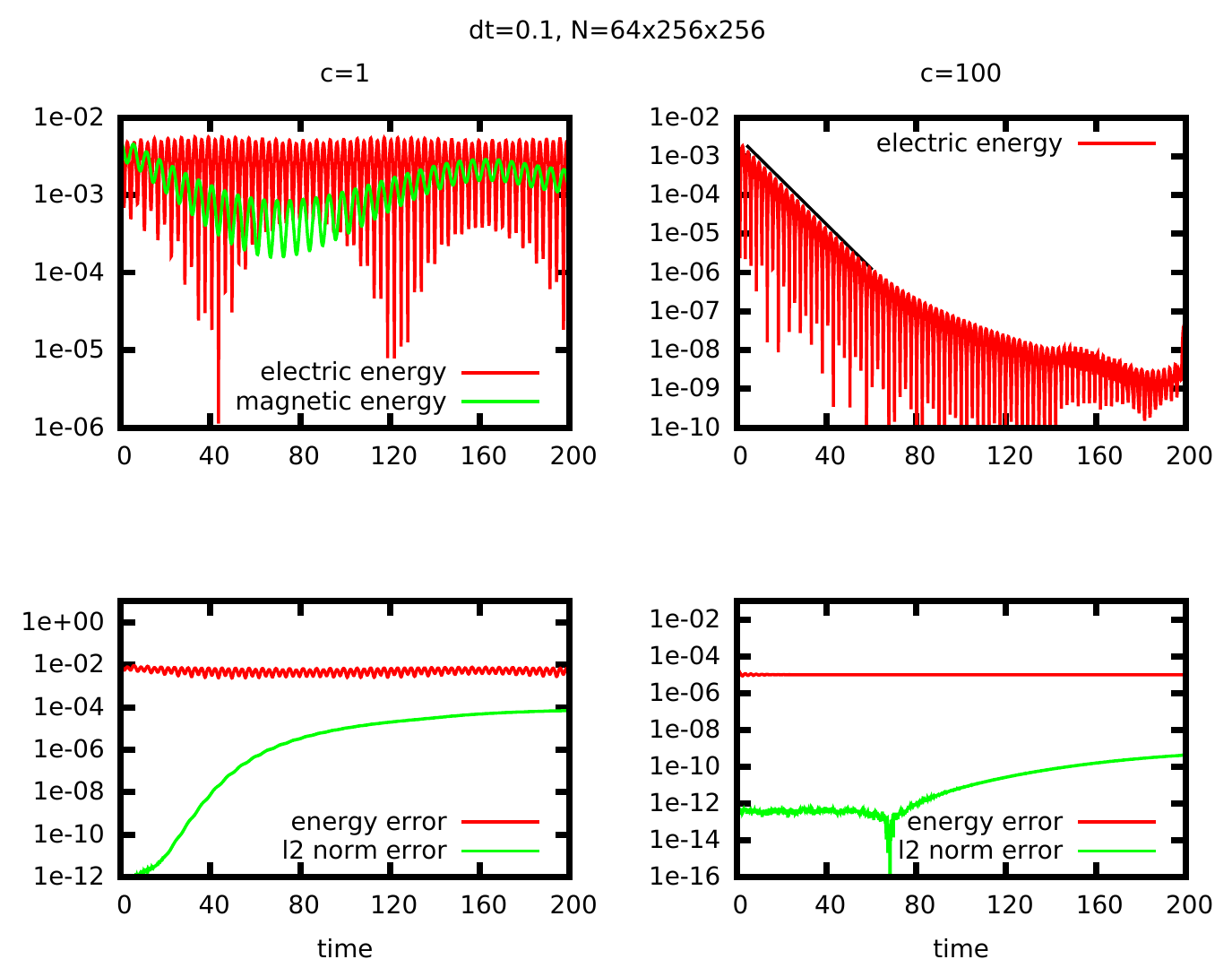}
	\caption{Fully relativistic case (Landau problem). 
	Top: time evolution of the  electric and magnetic energy obtained by 
AP-VM (for $c=1, 10^2$). 
Bottom: time evolution of the energy error and $L^2$ norm  (for $c=1, 10^2$). 
The analytic decay rate $(-0.0661)$ for the Vlasov--Poisson Landau damping is shown as a black line.
\label{fig:ldmagrel-evol}}
\end{figure}

We consider the same initial condition as in the semi-relativistic case, but now we set  $\gamma=\sqrt{1+|p|^2/c^2}$. The numerical parameters are $\Delta t=0.1, N_{x_1}=64, N_{p_1}=N_{p_2}=256$. 
As in the semi-relativistic case, we are interested in the time evolution of the electric and magnetic energies, 
in the error on the energy $\Hc(t)-\Hc(0)$ and in the error $\|f(t)-f(0)\|_{L^2}$ (in the $L^2$ norm in $x$ and $p$), 
for different values of $c$ ($c=1, 100$). The results are shown in Figure \ref{fig:ldmagrel-evol}.   
For the case $c=1$ we expect a complex interplay between the electric and magnetic field modes as well as with the plasma system. As we increase the dimensionless parameter $c$, however, Landau damping eventually dominates the dynamic of the system. For $c=100$ we in fact observe excellent agreement with the analytical decay rate that has been derived for the Vlasov--Poisson equations  
(see \cite{sonnen}).  This shows that our scheme converges to the correct limit in this example. In addition, we observe that the error in the total energy as well as the error in the $L^2$ norm decreases as we increase $c$. This might be due to the fact that the 2-dimensional interpolation in the $p$ direction 
degenerates as $c$ becomes large, so that $f$ remains unchanged during this step and does 
not affect the $L^2$ norm.

\subsubsection{Weibel type problem}
\begin{figure}
	\centering
	\includegraphics[width=12cm]{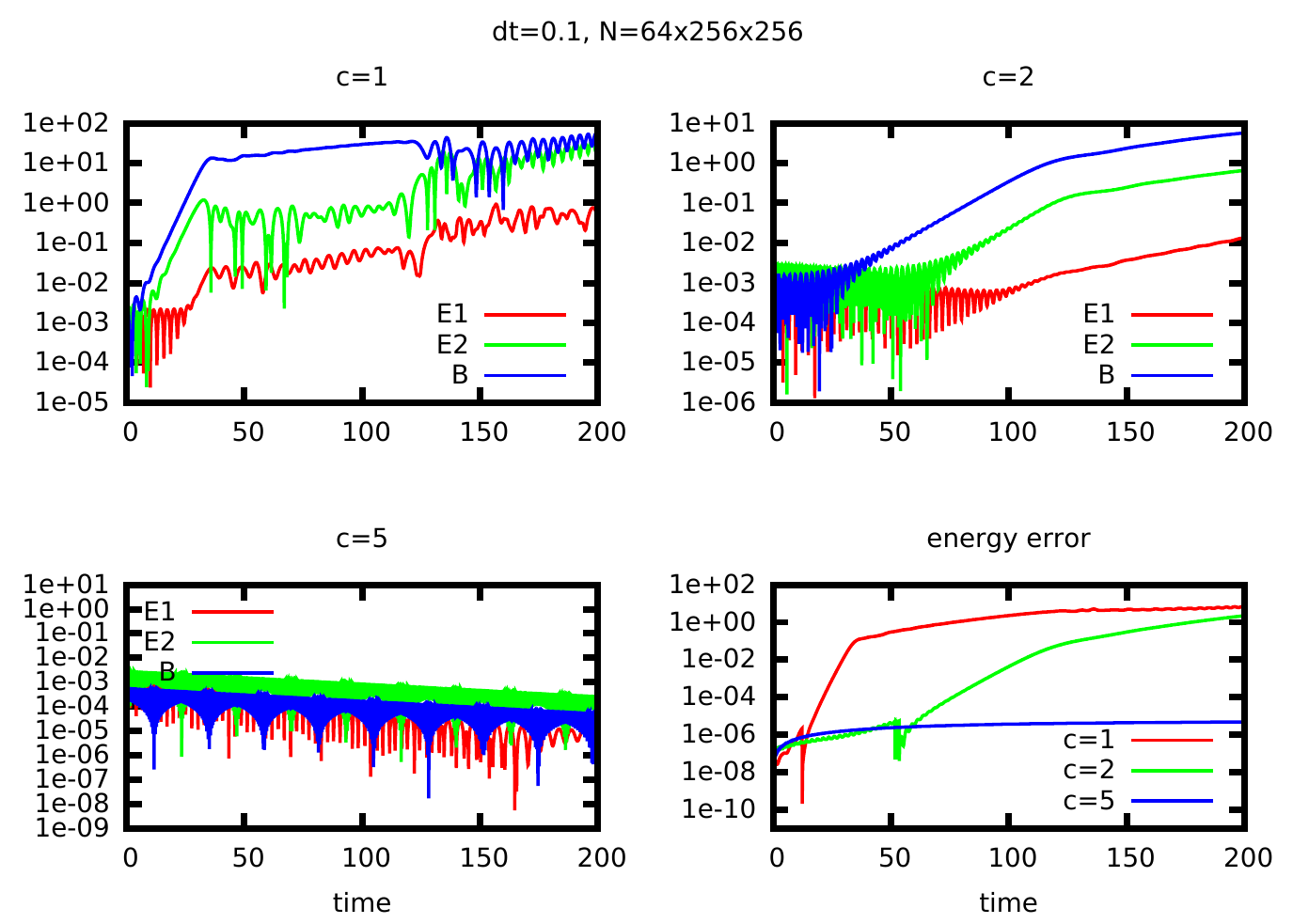}
\caption{Fully relativistic case (Weibel problem): time evolution of the most unstable mode ($k=1.25$) 
of the magnetic and the two electric fields obtained by AP-VM (for $c=1, 2, 5$). Bottom right: time 
evolution of the energy error for $c=1, 2, 5$. \label{eq:ttrel-evol}}
\end{figure}

Let us also consider the Weibel instability for the fully relativistic case (i.e.,~where $\gamma=\sqrt{1+|p|^2/c^2}$). 
The same initial condition and diagnostics as in the semi-relativistic case are considered. 
The numerical results are shown in Figure \ref{eq:ttrel-evol}. The dynamic is distinct in the sense that we also observe a significant growth in the electric field mode, which makes this test more challenging. 
As for the $\gamma=1$ case the Weibel instability eventually ceases to exist as we increase the dimensionless parameter $c$. 

\section{Conclusion \label{sec:conclusion}}
In the present work, we did propose a new time integrator for the Vlasov--Maxwell system that is asymptotic preserving in the classical limit (i.e.,~when the Vlasov--Maxwell system degenerates to the Vlasov--Poisson system). 
The method is based on a splitting scheme for the Vlasov equation and an implicit integrator for the linear part of Maxwell's equations. The choice of the latter is in fact crucial in order to obtain a stable numerical scheme in the relevant limit regime (i.e.,~for large values of the dimensionless parameter $c$). %

Numerical simulations show that the asymptotic preserving scheme proposed in this paper can be applied without severe time steps restrictions even for very large values of $c$. This gives the scheme a decisive advantage in the relevant regime compared to traditional time integrators. We have conducted a number of simulations illustrating the correct limit behavior in the classical regime.  In addition, we have demonstrated that for $c=1$ the numerical scheme agrees with the analytically derived growth rate for the Weibel instability for sufficiently small time step sizes.

In summary, we have constructed a time integrator that combines the computational advantages of the %
splitting scheme derived in \cite{cef} with the asymptotic preserving property for the classical limit. Such a scheme is of interest for numerical simulations in which magnetic effects are relatively weak but where the dynamic goes beyond what can be simulated using the more commonly employed Vlasov--Poisson model.

\clearpage\section*{Appendix: Dispersion relation and linear analysis \label{sec:dispersion}}
In this section we derive the dispersion relation for the Vlasov--Maxwell equations  both for 
the continuous and semi-discrete case (discrete in time but continuous in space) 
in the semi-relativistic configuration ($\gamma=1$). 
The dispersion relation does rely on linear analysis and thus only captures phenomena 
which are close to a steady state solution. However, as they give an indication
on the stability of a given mode,
it is instructive to compare the dispersion relation for the exact 
solution with the one obtained for the asymptotic preserving scheme proposed in this paper. It should 
be emphasized that the linear analysis we are going to conduct has been extensively 
used in the physics literature (in the continuous case) in order to determine a variety 
of properties of the Vlasov--Maxwell and Vlasov--Poisson systems (see for instance \cite{cpbm, sonnen, suzuki}). %

\subsection{Continuous dispersion relation}
We linearize the Vlasov--Maxwell system around a steady state given by $f_0(p), E_1=0, E_2=0, B=0$. 
For example, the well known Maxwell--Boltzmann distribution fits into this framework 
as does the temperature anisotropic initial value \eqref{init_weibel}
considered for the Weibel instability. 
Introducing the first order perturbations 
$\delta f$, $\delta E$, and $\delta B$ the linearized Vlasov equation can be written as
$$
\partial_t \delta f + p_x\partial_x \delta f + \delta E_1 \partial_{p_x}f_0+ \delta E_2 \partial_{p_y}f_0 
-  p_x\delta B \partial_{p_y}f_0 + p_y\delta B \partial_{p_x}f_0 = 0.  
$$
We now perform the Fourier transform  of the Vlasov--Maxwell equations in the spatial variable $x$ and the Laplace transform  in time. 
For Maxwell's equations we obtain
\begin{equation}
\label{rel-continu-EB}
-i\omega \delta E_2 = -c^2 ik \delta B - \delta J_2, \;\; -i\omega \delta B = -ik \delta E_2, \;\; -i\omega \delta E_1 = -\delta J_1. 
\end{equation}
The Vlasov equation becomes
\begin{equation}
\label{lin-EB}
-i\omega \delta f + p_x ik \delta f + \delta E_1 \partial_{p_x}f_0+ \delta E_2 \partial_{p_y}f_0 
-  p_x \delta B \partial_{p_y}f_0 + p_y \delta B \partial_{p_x}f_0=0.  
\end{equation}
Using the relation $\delta B=(k/\omega) \delta E_2$, we get 
$$
-i\omega \delta f +  p_x ik \delta f + \delta E_1 \partial_{p_x}f_0+ \delta E_2 \partial_{p_y}f_0 
-  p_x (k/\omega) \delta E_2 \partial_{p_y}f_0 + p_y(k/\omega) \delta E_2 \partial_{p_x}f_0=0. 
$$
Neglecting the $\delta E_1$ term and grouping the remaining terms we get
\[ i(-\omega  +  p_x k )\delta f = \delta E_2\left[ \partial_{p_y}f_0\frac{-\omega+p_x k }{\omega }-\frac{p_y k }{\omega}\partial_{p_x}f_0\right] \]
which yields after some manipulation
\begin{eqnarray*}
\delta f&=& -\frac{i\delta E_2}{\omega}\left[  \frac{p_y k}{ \omega - p_x k} \partial_{p_x}f_0+\partial_{p_y}f_0\right]
\end{eqnarray*}
Inserting the above expression for $\delta f$ into Maxwell's equations (using $\delta J_2 = \int p_y \delta f dp$), 
we obtain 
\begin{equation}
\label{rel-continu-JE}
\begin{split}
\delta J_2 &= -\frac{i\delta E_2}{\omega}\left[ \int \frac{p_y^2 k \partial_{p_x} f_0}{ ( \omega - p_x k)} dp_x dp_y + \int  p_y \partial_{p_y} f_0 dp_x dp_y\right] \\
&=:- \frac{i\delta E_2}{\omega}\left[ L_1(k, \omega, c) + L_2(c)\right].
\end{split}
\end{equation}
Using $-i\omega \delta B = -ik \delta E_2$, we deduce from Amp\`ere's equation $-i\omega \delta E_2 = -c^2 ik \delta B - \delta J_2$ the following relation
$$
i(\omega -c^2 k^2/\omega) \delta E_2 =  \delta J_2, 
$$
which immediately gives the dispersion relation 
\begin{equation}
\label{rel_disp_1}
0=-\omega^2 + k^2 c^2 - L_1(k, \omega, c) -L_2(c). 
\end{equation}
Note that a relation between $\delta J_2, \delta E_2$ and $\delta B$ can be derived 
by integrating \eqref{lin-EB} with respect to $p$ (after multiplying by $p_y$). This yields
\begin{equation}
\label{rel-continu-JEB}
\delta J_2 = \frac{i\delta E_2}{k} L_1+\frac{i\delta B}{k} L_2, 
\end{equation}
where $L_1$ and $L_2$ are given by 
\begin{equation}
\label{L1L2}
L_1 = \int_{\mathbb{R}^2} p_y\frac{\partial_{p_y} f_0}{p_x-\omega/k}\dd p, \;\; L_2= \int_{\mathbb{R}^2} \frac{p_y^2\partial_{p_x} f_0 - p_x p_y \partial_{p_y} f_0}{p_x-\omega/k}\dd p. 
\end{equation}

We now consider  the initial value of the Weibel instability 
$$
f_0(p)=\frac{1}{\pi v_{th}^2 \sqrt{T_r}}\exp\left( -\frac{(p_x-a)^2}{v_{th}^2} -\frac{(p_y-b)^2}{v_{th}^2 T_r}   \right)
$$
to compute $L_1$ 
\begin{eqnarray*}
L_1 &=& \int \frac{p^2_y}{\sqrt{\pi T_r} v_{th}} \exp\left( -\frac{(p_y-b)^2}{v_{th}^2 T_r} \right)dp_y \\
&& \times  \int \frac{(p_x-a)}{\sqrt{\pi } v^3_{th} (p_x-\omega/k)} \exp\left( -\frac{(p_x-a)^2}{v_{th}^2 } \right)dp_x\nonumber\\
&=& \left(\frac{v^2_{th}T_r }{2} + b^2\right) \times I\nonumber\\
\end{eqnarray*}
where $I$ is given by 
\begin{eqnarray*}
I &=&  \int \frac{(p_x-a)}{\sqrt{\pi } v^3_{th} (p_x-\omega/k)} \exp\left( -\frac{(p_x-a)^2}{v_{th}^2 } \right)dp_x\nonumber\\
&=& \frac{2}{v_{th}^2}\left[1 + \frac{(\omega/k-a)}{v_{th}} Z\left( \frac{\omega/k-a}{v_{th}} \right) \right]\nonumber\\
\end{eqnarray*}
with $Z(\xi) = 1/\sqrt{\pi}\; \int_0^\xi e^{-u^2}du=\sqrt{\pi} \exp(-\xi^2)(i-\mbox{erfi}(\xi))$. We thus obtain for $L_1$
\begin{multline*}
 L_1(k, \omega, c) = \left(T_r+\frac{2b^2}{v_{th}^2}\right)\\
\times  \left[ 1+\frac{\omega/k-a}{v_{th}}\sqrt{\pi}\exp\left(-\frac{(\omega/k-a)^2}{v^2_{th}}  \right) \left(i-\mbox{erfi}\left(\frac{\omega/k-a}{v_{th}}\right)\right)  \right]. 
\end{multline*}
A simple calculation shows that $L_2=-1$. Thus, from \eqref{rel_disp_1}, 
the dispersion relation can be written as follows 
\begin{multline}
-1+\omega^2- k^2 c^2 \\Ê
+ \left(T_r+\frac{2b^2}{v_{th}^2}\right)\left[ 1+\frac{\omega/k-a}{v_{th}}\sqrt{\pi}\exp\left(-\frac{(\omega/k-a)^2}{v^2_{th}}  \right) \left(i-\mbox{erfi}\left(\frac{\omega/k-a}{v_{th}}\right)\right)  \right]\\
 =0, 
\end{multline}
In the following we consider the case where $a=0$ and $b=0$. Thus, determining the zeros of
\begin{multline*}
D(\omega, k):=\\-1+\omega^2-k^2 c^2+ T_r\left[ 1+\frac{\omega/k}{v_{th}}\sqrt{\pi}\exp\left(-\frac{(\omega/k)^2}{v^2_{th}}  \right) \left(i-\mbox{erfi}\left(\frac{\omega/k}{v_{th}}\right)\right)  \right]
\end{multline*}
for a fixed $k$, $v_{th}$, $T_r$ and $c$, allows us to determine the stable and unstable perturbation. More precisely, a $\omega$ with a negative imaginary part corresponds to an unstable mode the amplitude of which grows exponentially in time (at least in the regime of validity of the linear analysis).

\subsection{Semi-discrete dispersion relation}
We repeat the linear analysis of the previous section for our time discretization. For the first step of the splitting, we consider for simplicity the following explicit Euler scheme  
\begin{eqnarray*}
f^{\star} &=& f^{n} (1-ikp_x \Delta t), \nonumber\\
E_1^{\star}    &=&  \hat{E}_1^n -  \Delta t\int_{\R^2} p_x f^n  \dd p, \nonumber\\
E_2^{\star}    &=&  \hat{E}_2^n - \Delta t\int_{\R^2} p_y f^n  \dd p.    
\end{eqnarray*}
The second step, in the case of the the implicit Euler scheme, is given by
$$
E_2^{n+1}=\frac{1}{1+\Delta t^2c^2k^2}(E_2^\star-\Delta tc^2 ik B^n), \;\; B^{n+1}=\frac{1}{1+\Delta t^2c^2k^2}(-\Delta t i k E_2^\star+B^n). 
$$
The third and fourth steps are given by the solution of the following two equations 
$$
\partial_t f + E^{n+1}\cdot \partial_p f = 0, \;\;\;\partial_t f + p_y B^{n+1} \partial_{p_x} f -p_x B^{n+1} \partial_{p_y} f = 0. 
$$
Since these two steps are nonlinear, we consider in this linear analysis, the corresponding linearization
$$
\partial_t f + E_2^{n+1} \partial_{p_y} f_0 = 0, \;\;\;\partial_t f + p_y B^{n+1} \partial_{p_x} f_0 -p_x B^{n+1} \partial_{p_y} f_0 = 0, 
$$
which can be solved exactly
$$
f^{n+1} = f^{\star} - \Delta t E_2^{n+1} \partial_{p_y} f_0- \Delta t p_y B^{n+1} \partial_{p_x} f_0 +  \Delta t p_x B^{n+1} \partial_{p_y} f_0.
$$
Now, for $g^n=(f^n, E_1^n, E_2^n, B^n)$ we consider the following Ansatz
$$
g^n = \delta g \exp(-i\omega n\Delta t).  
$$
Then the first step of the splitting becomes
\begin{eqnarray*}
\delta f^\star &=& \delta f e^{-i\omega n \Delta t} (1-\Delta t ikp_x), \nonumber\\
\delta E_2^\star &=& \delta E_2e^{-i\omega n\Delta t}  -\Delta t \int_{\mathbb{R}^2} p_y \delta f e^{-i\omega n \Delta t}  \dd p.
\end{eqnarray*}
The second step becomes 
\begin{eqnarray*}
\delta E_2 e^{-i\omega \Delta t} &=& \frac{1}{1+\Delta t^2c^2k^2}(E_2^\star e^{i\omega n \Delta t} - \Delta tc^2 ik \delta B)\nonumber\\
&=& \frac{1}{1+\Delta t^2c^2k^2}\left[ \delta E_2 - \Delta t \int_{\mathbb{R}^2} p_y \delta f \dd p - \Delta tc^2 ik \delta B\right], \nonumber\\
\delta B e^{-i\omega \Delta t} &=& \frac{1}{1+\Delta t^2c^2k^2} (-\Delta t ik E_2^\star e^{i\omega n \Delta t}+ \delta B)\nonumber\\
&=& \frac{1}{1+\Delta t^2c^2k^2}\left[-\Delta t ik \delta E_2 + \Delta t^2 ik \int_{\mathbb{R}^2} p_y \delta f \dd p + \delta B\right].
\end{eqnarray*}
The final step becomes
\begin{eqnarray*}
\delta f e^{-i\omega \Delta t} &=& f^\star e^{i\omega n\Delta t} - \Delta t \delta E_2  e^{-i\omega \Delta t}\partial_{p_y}f_0 - \Delta t \delta B e^{-i\omega \Delta t} \left[p_y\partial_{p_x} f_0 -p_x\partial_{p_y} f_0\right],\nonumber
\end{eqnarray*}
so that using  $\delta f^\star = \delta f e^{-i\omega n \Delta t} (1-\Delta t ikp_x)$ we obtain
$$
\delta f(e^{-i\omega \Delta t} -1+\Delta t ikp_x)
=- \Delta t \delta E_2  e^{-i\omega \Delta t}\partial_{p_y}f_0 - \Delta t \delta B e^{-i\omega \Delta t} \left[p_y\partial_{p_x} f_0 -p_x\partial_{p_y} f_0\right], 
$$
which after some manipulation yields
$$
\delta f = -\Delta t\delta E_2 \frac{\partial_{p_y}f_0}{1-e^{i\omega \Delta t}(1-\Delta tikp_x)} -\Delta t\delta B \frac{p_y\partial_{p_x}f_0 - p_x\partial_{p_y}f_0}{1-e^{i\omega \Delta t}(1-\Delta tikp_x)}. 
$$
As before we are able to express the current as a function of the electric and magnetic field perturbations
\begin{eqnarray*}
\int_{\mathbb{R}^2} p_y \delta f \dd p &=& -\Delta t\delta E_2 \int_{\mathbb{R}^2} p_y \frac{\partial_{p_y}f_0}{1-e^{i\omega \Delta t}(1-\Delta tikp_x)} \dd p \\Ê&& -\Delta t\delta B \int_{\mathbb{R}^2} p_y\frac{p_y\partial_{p_x}f_0 - p_x\partial_{p_y}f_0}{1-e^{i\omega \Delta t}(1-\Delta tikp_x)}\dd p, \nonumber\\
&=&  -\frac{\Delta t\delta E_2}{\Delta tike^{i\omega \Delta t}} \int_{\mathbb{R}^2}p_y \frac{\partial_{p_y}f_0}{p_x-\frac{(1-e^{-i\omega \Delta t})}{\Delta t ik}} \dd p \\ && -\frac{\Delta t\delta B}{\Delta tike^{i\omega \Delta t}} \int_{\mathbb{R}^2} p_y \frac{(p_y\partial_{p_x}f_0 - p_x\partial_{p_y}f_0)}{p_x-\frac{(1-e^{-i\omega \Delta t})}{\Delta t ik}} \dd p \nonumber\\
&:=&  \frac{i\delta E_2 e^{-i\omega \Delta t}}{k} L_1^{\Delta t} +\frac{i\delta B e^{-i\omega \Delta t}}{k} L_2^{\Delta t}. 
\end{eqnarray*}

Hence, we obtain a $3$x$3$ linear system $A^{\Delta t} U=0$ with $U=(\delta J_2, \delta E_2, \delta B)$, where $\delta J_2= \int_{\mathbb{R}^2} p_y \delta f \dd p$ and
$$
A^{\Delta t} = 
\left(
\begin{array}{ccc}
1& - \frac{i}{k} e^{-i\omega \Delta t} L_1^{\Delta t}& -\frac{i}{k}e^{-i\omega \Delta t} L_2^{\Delta t}\\[2ex]
1 & \frac{1}{\Delta t}\left[e^{-i\omega\Delta t} - \frac{1}{1+\Delta t^2c^2k^2}\right] & c^2 ik\\[2ex]
\Delta t^2 i k & ik &  \frac{1}{\Delta t }\left[e^{-i\omega\Delta t} - \frac{1}{1+\Delta t^2c^2k^2}\right]
\end{array}
\right).
$$
The dispersion relation in the semi-discrete case is hence given by det$(A^{\Delta t})=0$.

At the continuous level, we can, using \eqref{rel-continu-EB} and \eqref{rel-continu-JEB}, write the dispersion relation in matrix form. The dispersion relation is given by the zeros of the determinant of the following matrix
$$
A=
\left(
\begin{array}{ccc}
1& -\frac{i}{k} L_1& -\frac{i}{k} L_2\\[2ex]
1 & -i\omega & c^2 i k\\[2ex]
0 & ik & -i\omega
\end{array}
\right), 
$$
where $L_1$ and $L_2$ are given by \eqref{L1L2}.
It is easy to verify that $A^{\Delta t} \rightarrow A$ as $\Delta t\rightarrow 0$. This shows that the semi-discrete dispersion relation converges to the continuous dispersion.

\subsection{Dispersion relation for the Weibel instability}
\label{appendix3}
\begin{figure}
	\centering
	\includegraphics[width=12cm]{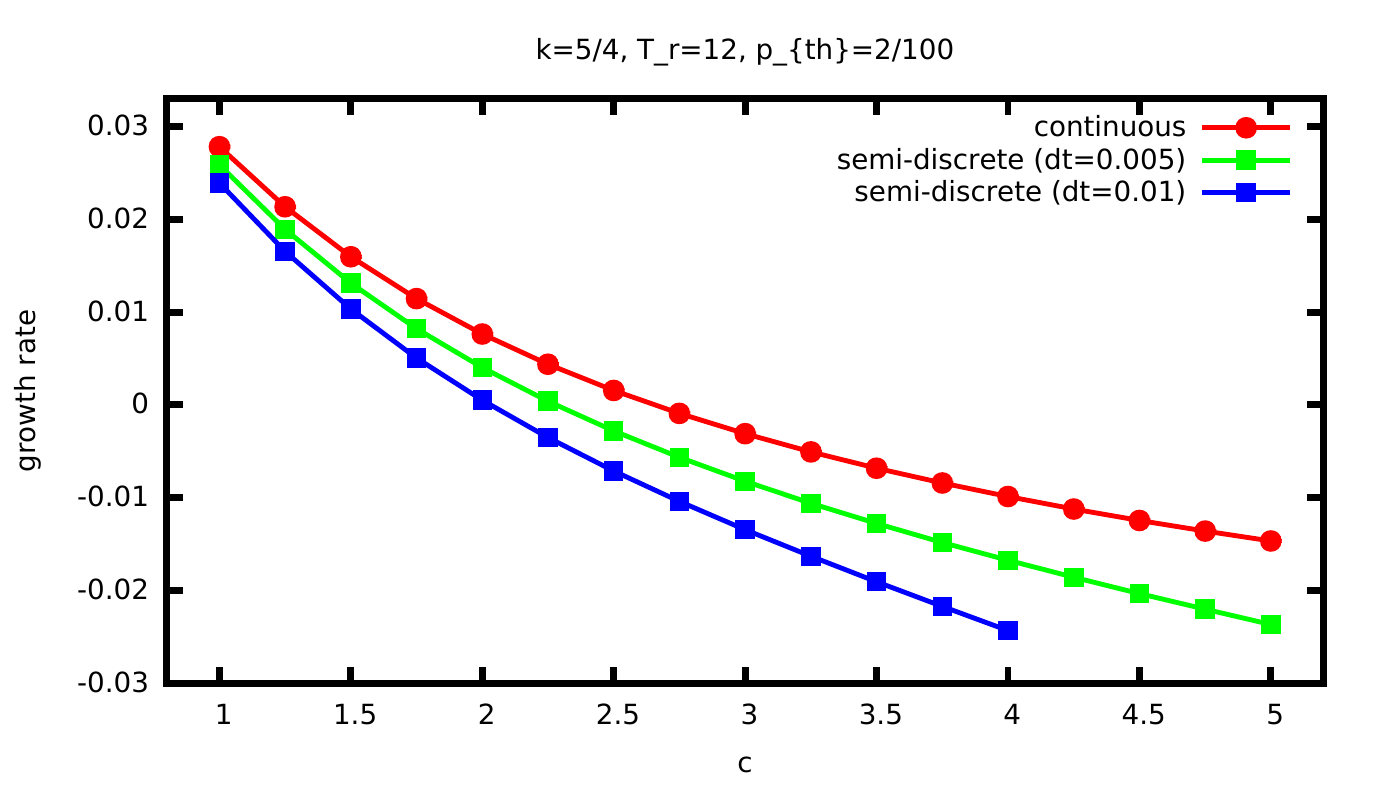}
	\caption{The growth rate of the Weibel instability for the continuous problem and the semi-discrete problem (for two different values of the step size $\tau$) are shown. \label{fig:dispersion}}
\end{figure}

The dispersion relations that have been derived for the continuous and the semi-discrete case are not amendable to a closed form solution. They can, however, be solved using a numerical root finding algorithm. The results for the parameters that have been used in the numerical simulation of the Weibel instability conducted in section \ref{sec:numerical-results} are shown in Figure \ref{fig:dispersion}.

We note that in order to obtain good agreement with the continuous formulation a relatively small time step size has to be chosen. This is something we already observed in the numerical simulations that have been conducted in section \ref{sec:numerical-results}. Let us also remark that for a value of $c$ above approximately $3$ the Weibel instability ceases to exist. In this regime the linear theory predicts a decay of the corresponding mode, which is also observed in the numerical simulations. 

\end{document}